\documentclass{article}

\usepackage{graphicx}      
\usepackage{amsmath,amssymb,mathrsfs}
\usepackage{enumerate}
\usepackage{bm}
\newtheorem{Theorem}{Theorem}
\newtheorem{Remark}{Remark}
\newtheorem{Lemma}{Lemma}
\newtheorem{Corollary}{Corollary}

\begin{document}
\title{Extended PID Control of Nonlinear Uncertain Systems}

\author{Cheng~Zhao and Lei~Guo\footnote{The authors are with the Key Laboratory of Systems and Control,
Academy of Mathematics and Systems Science, Chinese Academy of
Sciences, Beijing 100190, China, and Cheng~Zhao is also with School of Mathematical
Science, University of Chinese Academy of Sciences, Beijing 100049,
China (e-mail: zhaocheng@amss.ac.cn; Lguo@amss.ac.cn). Corresponding author: Lei Guo.},~\emph{Fellow,~IEEE}}
\date{}
\maketitle

\begin{abstract}
Since the classical proportional-integral-derivative (PID) controller is the most widely and successfully used ones in  industrial processes, it is of vital importance to investigate theoretically the rationale of this ubiquitous controller in dealing with nonlinearity and uncertainty. Recently, we have investigated the capability of the classical PID control for  second order nonlinear uncertain systems and provided  some analytic design methods for the choices of PID parameters, where the system is   assumed to be in the canonical form of cascade integrators. In this paper, we will consider the natural extension of the classical   PID control  for high order affine-nonlinear uncertain systems. In contrast to most of the literature on controller design of nonlinear systems, we do not require such special system structures  as  pure-feedback form, thanks to the strong robustness of the extend PID controller. To be specific,  we  will show that under some suitable  conditions on  nonlinearity  and uncertainty of the systems,    the extended  PID controller can globally(or semi-globally) stabilize the nonlinear uncertain systems, and at the same time the regulation error converges to $0$ exponentially fast,  as long as the control parameters are chosen from an open unbounded parameter manifold.
\newline

\end{abstract}



\section{Introduction}
Over the past 60 years,  remarkable progresses in modern control theory  have been made, e.g., numerous advanced control techniques including optimal control, robust control, adaptive control, nonlinear control, intelligent control, etc have been introduced and investigated. However,   the classical PID (proportional-integral-derivative) controller (or its variations), which has nearly 100 years of history,   is still the most widely and successfully  used one in engineering systems by far( see e.g., \cite{Astrom2006,samad2017survey}), which exhibits its lasting vitality.

In fact, a recent survey \cite{samad2017survey} shows that the PID controller has much higher impact rating than the advanced control technologies and that we still have nothing that compares with PID. However, it has also been reported that most of the practical PID loops are poorly tuned, and there is strong evidence that PID controllers remain poorly understood \cite{o2006pi}. Therefore, as pointed out in \cite{Astrom1995}, better understanding of the PID control may considerably improve its widespread practice, and so contribute to better product quality. This is the primary motivation of our theoretical investigation of the PID controller.

As is well-known, the PID controller has been investigated extensively in the literature by various scientists and engineers, but most focus on linear systems (e.g., \cite{Astrom1995,ho2003pid,silva2007pid}), albeit almost all practical systems are nonlinear with uncertainties. Therefore, to justify the remarkable effectiveness of the PID controllers for real world systems, we have to face with nonlinear uncertain dynamical systems and to understand the rationale and capability of this controller.

Recently,  we have given a theoretical investigation for the convergence and design of  PID controller for a basic class of nonlinear uncertain systems (see \cite{zhao2017capability}, \cite{zhao2017pid} and (\cite{krstic2017applicability}). For example, in \cite{zhao2017capability} we have shown that for second order nonlinear uncertain dynamical systems, one can select the three PID parameters to globally stabilize the closed-loop systems and at the same time to make the output of the controlled system converge to any given setpoint, provided that  the nonlinear uncertain functions satisfy a Lipschitz condition.  Moreover, necessary and sufficient conditions for the selection of the PID parameters have also been discussed and provided in \cite{zhao2017pid}.
These results have demonstrated theoretically that the classical PID controller does indeed have large-scale robustness with respect to both the uncertain system structure and the selection of the  controller parameters. However, in the  work of \cite{zhao2017capability} and \cite{zhao2017pid}, we have only provided global convergence results for  second-order nonlinear uncertain systems of integral cascades and there is no uncertainty in the controller channel.

Actually, in the area of nonlinear control, extensive researches have been conducted on the controller design (e.g.,\cite{krstic1995nonlinear,isidori2013nonlinear,khalil2002nonlinear,huang2004nonlinear,karafyllis2011stability}). For examples, the feedback linearization method  by using full knowledge of the system nonlinear functions (e.g., \cite{cheng1985global}, \cite{isidori2013nonlinear}), the backstepping approach for pure feedback forms in \cite{krstic1995nonlinear}, the extremum seeking methods for nonlinear uncertain systems(see e.g., \cite{killingsworth2006pid,tan2006non}), and many other interesting design methods for certain triangular forms(see, e.g.,\cite{xue2018performance,jiang2015control,jiang2018small,jiang1997small,huang2004nonlinear}), as well as for feedforward nonlinear systems, see e.g., \cite{marconi2000robust}.

  In this paper, we will consider a general class of single input single output(SISO) affine nonlinear uncertain systems. We will introduce  a  generalized concept called the extended PID controller(which is the  high-order extension the classical PID controller). We will show that for a large class of $n$-dimensional SISO affine nonlinear uncertain systems,   an $(n+1)$-dimensional parameter manifold can  be constructed, from which the extended  PID controller parameters can be arbitrarily chosen to globally or semi-globally stabilize the  nonlinear uncertain systems  with the regulation error converging to $0$ exponentially, even if  the system may not be globally transformed into the canonical form by coordinate transformation, thanks to the strong robustness of the extended PID controller as will be demonstrated in the paper. Moreover, in the case where the derivatives of the regulation error are not available,  the same results will also be established by incorporating a high-gain differential observer.

The rest of the paper is organized as follows.
In  Section II,  we will introduce the problem formulation. The main results are presented in  section III.   Section IV contains the proofs of the main theorems. Section V will conclude the paper with some remarks.


\section{Problem formulation}
\subsection{Notations}
We first introduce some notations and definitions to be used throughout this paper:

Let $x$ be a vector in  the $n$-dimensional Euclidean space $\mathbb{R}^n$, $P$ be an $m\times n$ matrix, and  $x^T$, $P^T$ denotes the transpose of  $x$ and  $P$ respectively.

Also, let $\|x\|$ denote the Euclidean norm of  $x$, and $\|P\|$ denote  the operator norm of the matrix $P$ induced by the Euclidean norm, i.e.,   $\|P\|=\sup_{x\in\mathbb{R}^n, \|x\|=1} \|Px\|$, which is known to be the largest singular value of $P.$

Let $z(t)$ be a function of time $t$, then $\overset{.}{z}(t)$ denotes the time derivative of $z(t)$. For simplicity, we oftentimes omit the variable $t$ whenever there is no ambiguity  in  the sequel.

Let $\Phi:\mathbb{R}^n\to \mathbb{R}^{n}$. For $z\in\mathbb{R}^n$, we denote $\Phi^{-1}(z)\overset{\triangle}{=}\{x\in\mathbb{R}^n: \Phi(x)=z\}$.

A map $\Phi:\mathbb{R}^n\to \mathbb{R}^n$ is called injective if $\Phi(x)\neq \Phi(y)$ for any $x\neq y$  and  is called surjective if $\Phi^{-1}(z)\neq \emptyset$  for any $z\in \mathbb{R}^n$.
 Moreover,  a  map $\Phi$ is called a global diffeomorphism on $\mathbb{R}^n$ if it is both injective and  surjective, and both  $\Phi$ and its inverse mapping (also denoted by $\Phi^{-1}$ for simplicity) are continuously differentiable.

Consider the following single-input-single-output(SISO) affine nonlinear system,
\begin{equation}\label{1}
\begin{cases}
\overset{.}{x}=f(x)+g(x)u\\
y=h(x)
\end{cases}
\end{equation}
where $f,g:\mathbb{R}^n\to \mathbb{R}^{n}$, $h:\mathbb{R}^n\to \mathbb{R}$ are sufficiently smooth mappings.

 The mappings $f,g:\mathbb{R}^n\to \mathbb{R}^{n}$ are  called smooth vector fields on $\mathbb{R}^n$. Let the coordinates of $x$ be $x_i$ and the components of $f$ and $g$ be $f_i$ and $g_i$ respectively, $i=1,\cdots,n.$ Define $L_f h(x)\overset{\triangle}{=}\sum_{i=1}^n \frac{\partial h}{\partial x_i}(x)f_i(x)$, which is called the Lie derivative of $h$ along with the vector field $f$. Let us further denote  $L_gL_f h(x)\overset{\triangle}{=} \sum_{i=1}^n \frac{\partial L_f h}{\partial x_i} g_i(x)$, $L_f^k h(x)\overset{\triangle}{=} L_fL_f^{k-1} h(x)$, $k\geq 1$, with $L_f^0 h(x)\overset{\triangle}{=}h(x)$.
 Moreover, the Lie bracket of the two vector fields $f$ and $g$ at $x$ is defined by $[f,g](x)=\frac{\partial g}{\partial x}(x)f(x)-\frac{\partial f}{\partial x}(x)g(x)$, which is a new  vector field on $\mathbb{R}^n$.  We also need the notation for vector fields $ad_f^{k}g$, which is defined recursively by $ad_f^{k}g=[f,ad_f^{k-1}g]$, $k\geq 1$, with $ad_f^0g=g$.

%
%
%
\subsection{Background}
Let us consider the system  (\ref{1}), where  $f:\mathbb{R}^{n}\to \mathbb{R}^{n}$, $g:\mathbb{R}^{n}\to \mathbb{R}^{n}$ and $h:\mathbb{R}^{n}\to \mathbb{R}$ are sufficiently smooth    unknown vector fields or unknown functions.

Let $y^*\in\mathbb{R}$ be a given setpoint. Our control objective is to design a   controller $u(t)$ to globally (or semi-globally) stabilize the  system (\ref{1}) and to achieve  asymptotic regulation $\lim_{t\to \infty} y(t)=y^*$.

To start with, let us assume that the relative degree of (\ref{1}) is $n$ at some point $x_0\in \mathbb{R}^n$, i.e.(see, e.g., \cite{isidori2013nonlinear}),
 $$L_g L_f^k h(x)=0,0\le k\le  n-2; L_g L_f^{n-1} h(x_0)\neq 0$$ in a neighborhood  of $x_0$. Let us define \begin{equation}\label{2}
\Phi(x)=(h(x),L_{f}h(x),\cdots,L_{f}^{n-1}h(x))^T.
\end{equation}
 Then there exists a neighborhood $U_{x_0}$ of $x_0$, such that $\Phi$ is a diffeomorphism on $U_{x_0}$. Under the coordinates transformation $z\overset{\triangle}{=}(z_1,\cdots,z_n)^T=\Phi(x)$,
  we have
 \begin{equation}\label{3.0}
\begin{cases}
\overset{.}{z}_1=z_2\\
~~~~\vdots\\
\overset{.}{z}_n=a(z)+b(z)u
\end{cases}
\end{equation}
where $y(t)=z_1(t)$, and for $z\in \Phi(U_{x_0})$, \begin{equation} a(z)=L_{f}^{n}h(x)=L_{f}^{n}h(\Phi^{-1}(z)),\end{equation}\begin{equation}b(z) =L_g L_{f}^{n-1}h(x)=L_g L_{f}^{n-1}h(\Phi^{-1}(z)), \end{equation} see, e.g., \cite{isidori2013nonlinear}.

If  $f,g,h$ were known, then one can chose   $k_1,\cdots, k_n$ such that the following matrix  $$\begin{bmatrix}0&1&\cdots&0\\ \vdots&\vdots&\ddots&\vdots\\0&0&\cdots&1\\-k_1&-k_2&\cdots&-k_n\end{bmatrix}$$ is  Hurwitz. Then under the state feedback control law \begin{equation}\label{u} u(x)=\frac{-L_{f}^{n}h(x)+v}{L_g L_{f}^{n-1}h(x)},\end{equation} with
\begin{equation}\label{v} v=-k_1(z_1-y^*)-k_2 z_2-\cdots-k_n z_n\end{equation}
the  system (\ref{1}) could be linearized to
\begin{equation}\label{111}
\begin{cases}
\overset{.}{z}_1 &=z_2\\
&~\vdots \\
\overset{.}{z}_n&=-k_1(z_1-y^*)-k_2 z_2-\cdots-k_n z_n\\
y&=z_1
\end{cases}
\end{equation}
If $(y^*,0,\cdots, 0)\in \mathbb{R}^n$ is sufficiently close to $\Phi(x_0)$, then it is not difficult to show that the closed-loop system (\ref{1}),(\ref{u}) and (\ref{v}) will satisfy $\lim_{t\to\infty}y(t)=y^*$ whenever the  initial state  $x(0)$ lies in some neighborhood $U\subset U_{x_0}$.

As is well-known, there are some fundamental limitations of the above method: The  state feedback control law (\ref{u}) needs the exact  information about both the structure and states of (\ref{1}), which is unrealistic  in most practical situations. Moreover, even if we know the exact information of (\ref{1}), the robustness of such designed  controller is still a concern. Furthermore,  the established theoretical results are usually local and  it is hard to estimate the corresponding local region.


 Thus, a challenging fundamental  problems is: Can we establish a theory  on global(or semi-global) stabilization and regulation by designing an output feedback controller  without using the precise mathematical model information?

\subsection{Problem Formulation}
Let us reconsider the system (\ref{1}), where  $f:\mathbb{R}^{n}\to \mathbb{R}^{n}$, $g:\mathbb{R}^{n}\to \mathbb{R}^{n}$ are smooth vector fields and $h:\mathbb{R}^{n}\to \mathbb{R}$ is a smooth function, which are all assumed to be unknown.

It is well known that the structure of the classical  PID controller is simple and  its implementation  does
not need precise mathematical models.
 In fact, in our previous work \cite{zhao2017pid}, we have shown that for second order nonlinear uncertain systems,  the classical PID control $u(t)=k_{1} e(t)+k_{0} \int_{0}^{t} e(s)ds+k_{2} \frac{de(t)}{dt}$  has the ability to achieve global stabilization, where  $e(t)=y^*-y(t)$ is the regulation error. We  have also shown  in \cite{zhao2017capability} that the  classical PID controller cannot stabilize dynamical systems which are described by ordinary differential equations with order $\geq 3$ even for linear time invariant systems. Therefore, it inspires us to introduce the concept ``extended PID controller'', which is defined by
\begin{equation}\label{6}
u(t)=k_{1} e(t)+k_{0} \int_{0}^{t} e(s)ds+k_{2} \overset{.}{e}(t)+\cdots+k_{n} e^{(n-1)}(t)
\end{equation}
where  $e(t)=y^*-y(t)$ is the regulation error, $\overset{.}{e}(t),\cdots, e^{(n-1)}(t)$ are the time derivatives of $e(t)$ up to the $(n-1)^{\text{th}}$ order.

From the definition (\ref{6}), we know that  the extended PID controller is an output feedback whose design does not need the precise model of the plant (\ref{1}), the control variable $u(t)$  is simply a weighted linear combination of the proportional, integral and derivative terms of the system regulation error, where the  weighting parameters $(k_0,k_1,\cdots,k_n)$ are called extended PID parameters.

We will in this paper investigate the capability of the extended PID controller (\ref{6}).  To be specific, the control objective is to understand when and how  the extended PID controller can guarantee that  the output $y(t)$ converges to a given reference value $y^*$ globally  (or semi-globally) with an exponential
rate of convergence, under the dynamic uncertainty $f(\cdot)$, the   control channel uncertainty $g(\cdot)$ and the state-output uncertainty $h(\cdot)$.

In this paper, we will rigrously show that the extended PID controller does indeed have the abovementioned nice properties,  even if the systems may not be transformed into the normal form globally by the coordinate transformation.

\section{The main results}
\subsection{Assumptions}

First, we introduce some notations that will be used throughout this paper. Let $y^*$ be a setpoint. Denote
\begin{align}\label{zzz}
z^*\overset{\triangle}{=}(y^*,0,\cdots,0)^T \in\mathbb{R}^n
\end{align}
and denote \begin{align}\label{HH}H(x)\overset{\triangle}{=}(F(x),G(x)),\end{align} where
\begin{align}\label{FG}
F(x)\overset{\triangle}{=}L_{f}^{n}h(x),~~ G(x)\overset{\triangle}{=}L_g L_{f}^{n-1}h(x).
\end{align}


Assume that  the system (\ref{1}) has uniform relative degree $n$, i.e.,  
$$L_gL_f^{i}h(x)=0,i=0,\cdots,n-2; ~~G(x)\neq 0, \forall x\in\mathbb{R}^n.$$

Since $G(x)$ is continuous and $G(x)\neq 0,~ \forall x\in\mathbb{R}^n$, we know that  the sign of  $G(x)$ cannot change. Therefore we introduce the following  assumption.

\textbf{Assumption (A):} System (\ref{1}) has uniform relative degree $n$. Furthermore,  the sign of $G(\cdot)$ is known and $G(x)$ is uniformly bounded away from zero. Without  loss of generality, we assume that $G(x)\geq \underline{b}>0$ for any $x\in \mathbb{R}^n$.

\begin{Remark} From Assumption ($A$), we know that $J_\Phi(x)$ is invertible for any $x\in \mathbb{R}^n,$ where $J_\Phi(x)$ denotes the Jacobian matrix of $\Phi$, (see e.g., \cite{isidori2013nonlinear}).   Under  the new coordinates $z=\Phi(x)$,
 the system (\ref{1}) transforms into  the  canonical form of cascade integrators (\ref{3.0}) locally.
Since the Jacobian matrix $J_\Phi(x)$ is invertible for any $x\in \mathbb{R}^n$, we conclude that $\Phi:\mathbb{R}^n\to \mathbb{R}^n$ is a locally  injective map, but may not be a global  diffeomorphism. Therefore,  the system (\ref{1}) may not be globally transformed into the normal form (\ref{3.0})  in general, unless the $n$ vector fields $(-1)^{i-1}ad_{\tilde{f}}^{i-1}\tilde{g}(x),i=1,\cdots,n$  are complete, where $\tilde{f}(x)=f(x)-\frac{F(x)g(x)}{G(x)}$ and $\tilde{g}(x)=\frac{g(x)}{G(x)}$, (see \cite{isidori2013nonlinear}).
\end{Remark}

Let $\tau_1,\tau_2:[0,\infty)\to [0,\infty)$ be  two     increasing functions   with $\limsup_{\rho\to0}  \frac{\tau_2(\rho)}{\rho}< \infty$, which will be used to  describe the model uncertainty quantitatively.

Specifically, let us introduce the following assumption.

\textbf{Assumption (B):} The functions $\Phi$ and $H$ defined respectively by (\ref{2}) and (\ref{HH}) satisfy

(i) $\|\Phi(x)\|\le \tau_1(\|x\|), \|H(x)\|\le \tau_1(\|\Phi(x)\|)$ for any $x\in \mathbb{R}^n$,

(ii) There exists $x^*\in \Phi^{-1}(z^*)$ such that the ``gap'' of $H$ at $x^*$ is bounded by that of $\Phi$ in the sense that $$\|H(x)-H(x^*) \| \le \tau_2(\|\Phi(x)-\Phi(x^*)\|)$$ holds for any $x\in\mathbb{R}^n$.

\begin{Remark}If the system is of the following normal form,
\begin{equation}\label{3}
\begin{cases}
\overset{.}{x}_1&=x_2\\
&~\vdots\\
\overset{.}{x}_{n-1}&=x_n\\
\overset{.}{x}_n&=a(x_1,\cdots,x_n)+b(x_1,\cdots,x_n)u(t)\\
y(t)&=x_1(t)
\end{cases},
\end{equation}  where $a(\cdot)$ and $b(\cdot)$ are unknown   smooth functions on $\mathbb{R}^n$.  Then  for any setpoint $y^*\in\mathbb{R}$,  we can verify     Assumptions (A) and   (B) as long as  the functions  satisfies  the following simplified assumption:

\textbf{Assumption (B'):} $b(x)\geq \underline{b}>0$, $|a(x)|+|b(x)|+\left\|\nabla a(x)\right\|+\left\|\nabla b(x)\right\| \le \rho(\|x\|),~~ \forall x\in\mathbb{R}^n$, for some constant $\underline{b}>0$ and some increasing function $\rho(\cdot)$, where $\nabla a(x)$ denotes the gradient of $a(\cdot)$.
 \end{Remark}

 \begin{Remark}
 For the system (\ref{3}), one can show that under Assumption (B'), the extended PID controller (\ref{6})  has the ability  to semi-globally stabilize the system(see Corollary 2). However,  for the general affine-nonlinear uncertain system (\ref{1}), Assumptions (A) and (B) are not sufficient. In fact, we will give an example in Appendix B to illustrate that for any $R>0$ and for any choices of the extended PID controller parameters, there always exists initial states $\|x(0)\|\le R$, such that the the  solution of the closed-loop system will have finite escape time, even though   both Assumptions (A) and (B) are satisfied. Therefore, we need to introduce certain additional assumptions. It turns out that either of the following two conditions can ensure the solution of the closed-loop system exists in the entire time interval $[0,\infty)$,  under appropriate choice of the controller  parameters and initial states.
\end{Remark}

\textbf{Assumption (C):} The growth rate of the inverse  of the Jacobian matrix $J_{\Phi}(x)$ satisfies
    \begin{align}\label{13}
    \big\|J_\Phi ^{-1}(x)\big\| \le N_1\|x\|\log \|x\|+N_2,\forall x\in\mathbb{R}^n
    \end{align}   for some  constants $N_1>0, N_2>0$ (possibly unknown).

\textbf{Assumption (C'):} The coordinates transformation map $\Phi$ is a global diffeomorphism on $\mathbb{R}^n$.

\begin{Remark} We remark that for systems (\ref{3}), both Assumptions (C) and   (C') are satisfied since the coordinate transformation map $\Phi$ is the identity map, i.e., $\Phi(x)=x,~\forall x\in\mathbb{R}^n.$ On the other hand, we point out that  the  growth rate of  $\big\|J_\Phi ^{-1}(x)\big\|$ in Assumption (C) cannot be relaxed slightly  to, e.g.,
\begin{align}
\label{exa}\big\|J_\Phi ^{-1}(x)\big\| \le N_1\|x\|\log^{1+\eta}\|x\|+N_2, ~\forall x\in\mathbb{R}^n
\end{align} for any $\eta>0$.  We will give an example to illustrate this in  Appendix B.
\end{Remark}

\subsection{Semi-global stabilization}
We will first show that the extended PID controller  (\ref{6}) can stabilize the system $(\ref{1})$ semi-globally under the above assumptions in the following sense: For any $R>0$, there exists an open unbounded set $\Omega(R)\subset\mathbb{R}^{n+1}$, such that whenever the extended PID parameters $(k_0,\cdots,k_n)\in \Omega(R)$, then the closed-loop system will satisfy $\lim_{t\to \infty} y(t)=y^*$ for any initial state $x(0)$ satisfying $\|x(0)\|\le R$.
 To be specific, we have the following theorem.

 \begin{Theorem} Consider the SISO affine-nonlinear uncertain  system $(\ref{1})$ with the extended PID controller defined by (\ref{6}).  Suppose that Assumptions (A), (B) and (C) are satisfied, where $\underline{b}>0$ and  the increasing functions $\tau_1(\cdot),\tau_2(\cdot)$ are known.
Then for any $R>0$, there exists an $(n+1)$-dimensional unbounded parameter manifold $\Omega(R)\subset\mathbb{R}^{n+1}$,
 such that whenever $(k_0,\cdots,k_n)\in \Omega(R)$,  the solution of the the closed-loop system with initial state   $\|x(0)\| \le R$ will  exist in $[0,\infty)$ and the regulation error $e(t)$ will converge to zero exponentially. Moreover,    if Assumption (C) is replaced by Assumption (C'), then all the above statements still hold and   there exists $x^*\in \mathbb{R}^n$  depending  on $y^*$ and $\Phi$ only, such that the state $x(t)$ is bounded and $\lim_{t\to\infty}x(t)=x^*.$
\end{Theorem}

\begin{Remark}
The concrete construction of the parameter manifold $\Omega(R)$ can be found in the proof of Theorem 1 to be given in Section IV. Moreover, from the proof of Theorem 1, we can also get an upper bound of $e(t)$: whenever  $(k_0,\cdots,k_n)\in \Omega(R)$, we have for any $t\geq 0$,
\begin{align}\label{et}|e(t)|\le& \|\Phi(x(t))-z^*\|
 \le c e^{-\alpha t} \left\{\|\Phi(x(0))-z^*\|+\tau_1(|y^*|)/(k_0\underline{b})\right\},\end{align} for any initial states  satisfying $\|x(0)\| \le R$, where  $\alpha>0$ is a constant depends on $(k_0,\cdots,k_n)$ and $c>0$ is a constant only depends on $n.$ As a consequence, we have  $\lim_{t\to \infty}\Phi(x(t))=z^*$.
 Furthermore, if the initial state satisfies $\Phi(x(0))=z^*$, then from (\ref{et}), it is easy to see that $\sup_{t\geq 0} |e(t)|\le c\tau_1(|y^*|)/(k_0\underline{b}).$ Therefore, for any $\epsilon>0$, we can get $\sup_{t\geq 0} |e(t)|\le \epsilon$ as long as the integral parameter $k_0$ is suitably large.
\end{Remark}

 We remark that Assumption (B) is crucial for stabilization of (\ref{1}) by the extended PID controller, though it may not be easy to be verified in general. However, if the norm of the coordinate transformation  $\Phi(x)$ has a lower bound function and  $\|\Phi(x)\|,|F(x)|, |G(x)|, \|\nabla F(x)\|,\|\nabla G(x)\|$ have some upper bound functions, then the verification can be considerably simplified. To this end, we introduce the following  assumption on $\Phi$ and $H$ to replace assumption (B), which does not depend on the setpoint $y^*$.

\textbf{Assumption (B0):} The functions $\Phi$ and $H$ defined  by (\ref{2}) and (\ref{HH}) satisfy the following inequalities:
$$\|\Phi(x)\|+\|H(x)\|+\|J_\Phi^{-1}(x)\| +\|J_H(x) \| \le \rho_1(\|x\|)
,~\|\Phi(x)\|\geq \rho_0(\|x\|),~ ~~\forall x\in\mathbb{R}^n$$
where  $\rho_0(\cdot),\rho_1(\cdot)$ are two known continuous increasing  functions with $\lim_{r\to \infty}\rho_0(r)=\infty$.

By Theorem 1,  we can  get the following corollary which does not need Assumptions (C) or (C'):

 \begin{Corollary}
 Consider the SISO affine-nonlinear uncertain  system $(\ref{1})$ with the extended PID controller defined by (\ref{6}).  Suppose that Assumptions (A), (B0) are satisfied. Then for any setpoint $y^*\in \mathbb{R}$ and for any  $R>0$, there exists $x^*\in \mathbb{R}^n$ and an $(n+1)$-dimensional parameter manifold $\Omega\subset\mathbb{R}^{n+1}$, such that whenever $(k_0,\cdots,k_n)\in \Omega$,   the solution of the closed-loop system with initial state   $\|x(0)\| \le R$ will exist in $[0,\infty)$ and the system state $x(t)$ will be bounded with $\lim_{t\to\infty} x(t)=x^*$ exponentially, and at the same time the regulation error $e(t)$  converges to $0$ exponentially.
\end{Corollary}

We remark that the conditions used in the above corollary can be further simplified for  the basic class of  $n^{\text{th}}$-order uncertain chain of integrators (\ref{3}), since in this case the coordinate transformation map $\Phi(x)=x$.  The following corollary can be deduced  from Corollary 1 immediately.

\begin{Corollary} Consider the nonlinear uncertain system (\ref{3}) with the extended PID controller  (\ref{6}). Then for any setpoint $y^*$ and for any $R>0$, there exists an open unbounded set $\Omega\subset \mathbb{R}^{n+1}$, such that whenever $(k_0,\cdots,k_n)\in \Omega$,   the solution of the closed-loop system with initial state   $\|x(0)\| \le R$ will satisfy $\lim_{t\to\infty} e(t)=0$ and $\lim_{t\to\infty} x(t)=z^*$ exponentially, provided that the nonlinear uncertain  functions $a(\cdot)$ and $b(\cdot)$ satisfy Assumption (B').
\end{Corollary}

\subsection{Global stabilization}

In this section, we will show that the extended PID controller can globally stabilize the system $(\ref{1})$ under some  additional assumptions  on the mappings $\Phi,F$ and $G$.
Specifically,  if $G(x)$(which is defined in (\ref{FG})) has a constant upper bound $\overline{b}$ and the increasing function $\tau_2$ in Assumption (B) is a  linear function, i.e., $\tau_2(r)=L r$,  then the extended PID controller can globally stabilize  (\ref{1}) in the sense that:   there exists an open unbounded  parameter set $\Omega\subset\mathbb{R}^{n+1}$, such that the closed-loop system will satisfy $\lim_{t\to \infty } e(t)=0$ for any initial state $x(0)\in \mathbb{R}^n$ as long as $(k_0,\cdots,k_n)\in \Omega$.

Let $\tau:[0,\infty)\to [0,\infty)$ be an     increasing function and $L>0$, $\overline{b}>0$ be two constants. We introduce the following assumption to  describe the model uncertainty.

\textbf{Assumption (B1):} The functions $\Phi$ and $H$ defined respectively by (\ref{2}) and (\ref{HH}) satisfy:

 (i) $~G(x)\le \overline{b},$ $  \|H(x)\|\le \tau(\|\Phi(x)\|) \text{~for any~}  x\in \mathbb{R}^n,$

(ii) There exists $x^*\in \Phi^{-1}(z^*)$ such that$$\|H(x)- H(x^*)\| \le L \|\Phi(x)-\Phi(x^*)\|, ~\forall   x\in \mathbb{R}^n.$$
  Now, we have the following  global results on the extended PID controller:
\begin{Theorem}  Consider the SISO affine-nonlinear uncertain  system $(\ref{1})$ with the extended PID controller defined by (\ref{6}). Suppose that Assumptions (A), (B1) and (C) are satisfied. Then there exists an $(n+1)$-dimensional unbounded parameter manifold $\Omega\subset\mathbb{R}^{n+1}$, such that whenever $(k_0,\cdots,k_n)\in \Omega$,   the solution of the closed-loop system will exist in $[0,\infty)$ for any initial state $x(0)\in\mathbb{R}^n$ and that the regulation error $e(t)$ will converge to $0$ exponentially. Moreover,    if Assumption (C) is replaced by Assumption (C'), then all the above statements still hold and   there exists $x^*\in \mathbb{R}^n$  depending  on $y^*$ and $\Phi$ only, such that the state $x(t)$ is bounded and $\lim_{t\to\infty}x(t)=x^*.$
\end{Theorem}

 \begin{Remark} We first remark that, in contrast to the semi-global results established in Theorem 1, the global results in Theorem 2 do not require  the  parameters $(k_0,\cdots,k_n)$ of the extended PID   depends on the range of initial state $x(0)$. Second, we note that  from  the proof of Theorem 2, the second condition (ii) in Assumption (B1) can be slightly weakened to $|F(x)G(x^*)- F(x^*)G(x)| \le L \|\Phi(x)-\Phi(x^*)\|, ~\forall   x\in \mathbb{R}^n$, for some $L>0.$
\end{Remark}

Now, let us give three typical examples to show  how the conditions used in Theorem 2 can be further simplified if some additional information on the systems structure are available.

$\textbf{Example ~1 }(Nonlinear ~canonical~ form)$: Consider the nonlinear uncertain system (\ref{3}) with the extended PID controller (\ref{6}).
Then for any setpoint $y^*\in\mathbb{R}$, there exists an $(n+1)$-dimensional unbounded parameter manifold $\Omega\subset\mathbb{R}^{n+1}$, such that whenever $(k_0,\cdots,k_n)\in \Omega$,  the closed-loop system will satisfy
\begin{align*}\lim_{t\to \infty}e(t)=0,~ \lim_{t\to \infty}x(t)=z^*\end{align*}  for any initial state $x(0)\in \mathbb{R}^n$, provided that  the nonlinear uncertain functions $a(\cdot)$ and $b(\cdot)$ satisfy the following conditions respectively:
 $$~(\text{i}) ~|a(x)-a(y)|  \le L\|x-y\|, ~\forall   x,y \in \mathbb{R}^n~ \text{and} ~~|a(0)|\le M$$
  $$~(\text{ii})~ 0<\underline{b}\le b(x)\le \overline{b},~ |b(x)-b(y)|  \le L\|x-y\|, ~\forall x,y \in \mathbb{R}^n,$$
 where $\underline{b},\overline{b},L$ and $M$ are  known positive constants.

In the following two typical examples, we give some simple conditions for the classical PID to globally stabilize second order nonlinear uncertain systems, which generalize the  nonlinear models investigated in \cite{zhao2017pid}.

$\textbf{Example ~2 }(Second ~order ~system ~with~ strict-feedback~ form)$: The following  nonlinear uncertain system
\begin{equation}\label{triangular}
\begin{cases}
\overset{.}{x}_1=f_1(x_1)+x_2\\
\overset{.}{x}_2=f_2(x_1,x_2)+u(t)\\
y(t)=x_1(t)
\end{cases}
\end{equation}
can be globally stabilized by the classical PID controller,   provided that  the nonlinear uncertain functions $f_1$ and $f_2$ satisfies $|f_1'(x_1)|\le L,\forall x_1\in\mathbb{R}$ and $ |f_2(x)-f_2(y)| \le L\|x-y\|,\forall x,y\in\mathbb{R}^2.$ This example also generalizes the corresponding results in \cite{zhao2017capability}.

$\textbf{Example ~3 }(Second ~order ~system ~with~ pure-feedback~ form)$: The following $2$-dimensional  SISO system
\begin{equation}\label{triangularr}
\begin{cases}
\overset{.}{x}_1=f_1(x_1,x_2)\\
\overset{.}{x}_2=f_2(x_1,x_2)+g(x_1,x_2)u(t)\\
y(t)=x_1(t)
\end{cases}
\end{equation}
can be globally stabilized by the classical PID controller  as long as the following conditions are satisfied: $\left\|\nabla f_i(x)\right\|\le L$, $i=1,2$, $\left\|\nabla g(x)\right\|\le L$, $\|\text{Hess} (f_1)(x) \|\le L$
, $\frac{\partial f_1}{\partial x_2}(x)\geq \underline{b}, b_1\le  g(x)\le b_2$ for any $x=(x_1,x_2)\in\mathbb{R}^2$, and $ |f_1(0)|+|f_2(0)|\le M$,
 where   $\text{Hess} (f_1)$ is the Hessian matrix of $f_1$  and the constants $0<\underline{b}\le L$, $0<b_1\le b_2$ and $M>0$ are known.
\subsection{Extended PID controller with differential trackers}
From the definition of   the extended PID controller (\ref{6}), we know that the implementation of (\ref{6}) needs the derivatives    information  $\overset{.}{e}(t),\cdots, e^{(n-1)}(t)$ of the regulation error, if the system order $n\geq 2$. However, in most practical situations,  these derivatives  may not be available directly. Therefore, we need to construct a differential observer to obtain an online estimation of the derivatives  of the regulation error.

 In this subsection, we introduce the following high gain differential observers(see e.g., \cite{huang2014active},\cite{khalil2014high}):
\begin{equation}\label{5}
\begin{cases}
\overset{.}{\hat{z}}_1&=\hat{z}_2+\frac{\beta_1}{\epsilon}(e-\hat{z}_1)\\
&~\vdots\\
\overset{.}{\hat{z}}_{n-1}&=\hat{z}_n+\frac{\beta_{n-1}}{\epsilon^{n-1}}(e-\hat{z}_1)\\
\overset{.}{\hat{z}}_n&=\frac{\beta_n}{\epsilon^n}(e-\hat{z}_1)
\end{cases}
\end{equation}
where the parameters $(\beta_1,\cdots,\beta_n)\in\mathbb{R}^n$ are given  such that the polynomial $s^n+\beta_1 s^{n-1}+\cdots+\beta_n$ is Hurwitz and $\epsilon > 0$ is the observer gain parameter to be determined. We introduce the following  differential observer-based extended PID controller:
\begin{equation}\label{4}
u(t)=k_{1} e(t)+k_{0} \int_{0}^{t} e(s)ds+k_{2} \hat{z}_2(t)+\cdots+k_{n} \hat{z}_n(t)
\end{equation}
where  $\hat{z}_i(t),i=2,\cdots,n$ are the estimators of the derivatives $\frac{de(t)}{dt},\cdots, \frac{d^{n-1}e(t)}{dt^{n-1}}$ of the regulation error respectively.

Now, let us consider the closed-loop system $(\ref{1})$ with the extended PID controller defined by (\ref{5})-(\ref{4}).

\begin{Theorem} Consider the SISO affine-nonlinear uncertain  system $(\ref{1})$, where $u(t)$ is the differential observer-based extended PID controller.  Suppose that Assumptions (A), (B1) and (C) are satisfied. Then there exists an open unbounded set $\Omega\subset\mathbb{R}^{n+1}$, such that for any $(k_0,\cdots,k_n)\in \Omega$,  there exists $\epsilon^*>0$, such that for any $0<\epsilon <\epsilon^*$, and  for any initial states $x(0)\in\mathbb{R}^n$, $\hat{z}(0)\in \mathbb{R}^n$, the solution of the closed-loop system will   exist in $[0,\infty)$  and  the regulation error $e(t)$ will converge to zero exponentially.
 Moreover,    if Assumption (C) is replaced by Assumption (C'), then all the above statements still hold and   there exists $x^*\in \mathbb{R}^n$  depending  on $y^*$ and $\Phi$ only, such that the state $x(t)$ is bounded and $\lim_{t\to\infty}x(t)=x^*.$
\end{Theorem}

\begin{Remark}First, we remark that  from the proof of Theorem 3 to be given in the next section, a precise upper bound on the regulation error can be obtained. Second, we remark that semi-global results like Theorem 1  can also obtained, which will be discussed  in details elsewhere.\end{Remark}

\section{Proofs of the main results}
Before proving the theorems, we first list some lemmas, whose proofs are presented in Appendix A.

Denote $\lambda\overset{\triangle}{=}(\lambda_0, \cdots , \lambda_n)\in\mathbb{R}^{n+1}$ and  define an open unbounded set $\Omega_1\subset\mathbb{R}^{n+1}$ as follows:
 \begin{equation}\label{omega1}\Omega_1= \{\lambda \big | 2<\lambda_i-2i <3, i=0,\cdots,n-1; \lambda_n>2n+2 \}\end{equation}
 and for $\lambda\in \Omega_1$, we define a $(n+1)\times(n+1)$ matrix $P=P(\lambda)$ as follows (see \cite{cheng2004note}): \begin{equation}\label{P} P=\begin{bmatrix}(-\lambda_{0})^{-n}&\cdots&(-\lambda_{n})^{-n}\\
\vdots&&\vdots\\(-\lambda_0)^{-1}&\cdots&(-\lambda_{n})^{-1}\\1&\cdots&1\end{bmatrix}
\end{equation}
and denote $(d_0,\cdots,d_n)^T$ be the last column of $P^{-1}$, i.e., \begin{equation}\label{d} (d_0,\cdots,d_n)^T=P^{-1}(0,\cdots,0,1)^T
\end{equation}

\begin{Lemma} Under the above notations, let us define
\begin{align}\label{c1}
c_1&\overset{\triangle}{=}\sup_{\lambda\in \Omega_1}\|P\|,~~ c_2 \overset{\triangle}{=}\sup_{\lambda\in \Omega_1} \|P\|\|P^{-1}\|,
 c_3 \overset{\triangle}{=} \sup_{\lambda\in \Omega_1} \sqrt{n}(2n+1) d_n
 \end{align}
\begin{align}\label{c2}
c_4(i)\overset{\triangle}{=}\sup_{\lambda\in \Omega_1}\left|(2n+1) n \lambda_{n} d_i \right|, i=0,\cdots,n-1,
\end{align}
and denote \begin{align}\label{c0} c_0=\max\{c_1,c_2,c_3,c_4(i),i=0,\cdots,n-1\}, \end{align} then $c_0$ is a positive finite number.
\end{Lemma}

Note that Lemma 1 is nontrivial because $\lambda_n$ together with $\Omega_1$ are unbounded. To introduce other lemmas, we now define a parameter manifold first.
Let $c\geq c_0$  be any constant.
For $L>0$ and $0<\underline{b}\le\overline{b}$, we define the following $n+1$ dimensional parameter manifold $\Omega_{L,\underline{b},\overline{b},c}\subset\mathbb{R}^{n+1}$(which is open and unbounded in $\mathbb{R}^{n+1}$),
\begin{equation}\label{Omega}\Omega_{L,\underline{b},\overline{b},c}\overset{\triangle}{=}\left\{\begin{bmatrix} k_0\\ \vdots \\k_n \end{bmatrix} \bigg|  \begin{bmatrix} k_0\\ \vdots \\k_{n-1} \\k_n\end{bmatrix}=\frac{1}{\underline{b}}\begin{bmatrix}
\prod_{i=0}^{n}\lambda_{i}\\
\vdots\\\sum_{ i<j}\lambda_{i}\lambda_{j}\\
\sum_{i=0}^{n}\lambda_{i}
\end{bmatrix} , \lambda \in \Omega_{\Lambda} \right\}
\end{equation}
where  $\Omega_{\Lambda}$ is defined by
 \begin{align} \label{omega2}
 \Omega_{\Lambda}=\left \{\lambda\in \Omega_1 \bigg | \lambda_n> \left( L c^2+\frac{\left(\overline{b}-\underline{b}\right)c}{\underline{b}}\right)^2+Lc^2 \right\}\end{align}

In the following lemmas, the constant $T$  can be a finite positive number $0<T<\infty$   or an infinity $T=\infty.$

\begin{Lemma} Let $\overline{Y}(t)=(y_0(t),\cdots,y_n(t))^T$ be a continuously differentiable vector valued function on $[0,T)$. Suppose that  there exists $a_t$ and $b_t$ such that the following equalities hold for $ t \in [0,T)$,
\begin{equation}\label{close}
\begin{cases}
\overset{.}{y}_0&=y_1\\
&~\vdots\\\overset{.}{y}_{n-1}&=y_n\\
\overset{.}{y}_n&=a_t-b_t(k_0y_0+\cdots+k_n y_n)
\end{cases}
\end{equation}
where $|a_t|\le L\left\|\overline{Y}(t)\right\|$ and $0<\underline{b}\le b_t\le \overline{b}$, for any $t\in [0,T)$.  Then for any  $(k_0,\cdots,k_n)\in \Omega_{L,\underline{b}, \overline{b},c}$( where $\Omega_{L,\underline{b}, \overline{b},c}$  is defined in (\ref{Omega})), there exists $\alpha>0$, such that $\overline{Y}(t)$  satisfies
 $$\left\|\overline{Y}(t)\right\|\le  c e^{-\alpha t}\left \|\overline{Y}(0)\right\|, ~~\forall t\in [0,T)$$
\end{Lemma}
\begin{Lemma}
Consider the system of equalities (\ref{close}) again,
but where $|a_t|\le \tau_2(\left\|\overline{Y}(t)\right\|)$ and $0<\underline{b}\le b_t\le \tau_1(\left\|\overline{Y}(t)\right\|)$, for any $t\in [0,T)$ and where $\tau_1,\tau_2:\mathbb{R}^+\to \mathbb{R}^+$ are two     increasing functions with $\limsup_{\rho\to0}  \frac{\tau_2(\rho)}{\rho}< \infty$.  Then for any $R>0$, and any $(k_0,\cdots,k_n)\in \Omega_{L_0,\underline{b},b_0,c}$ with $L_0=\sup_{0\le \rho \le cR}\frac{\tau_2(\rho)}{\rho}$, $b_0=\tau_1(cR)$, there exists $\alpha>0$, such that  $\overline{Y}(t)$ satisfies
 $$\left\|\overline{Y}(t)\right\|\le  c e^{-\alpha t}\left \|\overline{Y}(0)\right\|,\forall t\in [0,T)$$
  provided  that $\left\|\overline{Y}(0)\right\|\le R$.
\end{Lemma}

\begin{Lemma} Let $\overline{Y}(t)=(y_0(t),\cdots,y_n(t))^T$, $\xi(t)=(\xi_1(t),\cdots,\xi_n(t))^T$ be two continuously differentiable vector valued functions on $[0,T)$. Suppose that $\forall t \in [0,T)$, there exists $a_t$ and $b_t$ such that the following equalities hold,
\begin{equation}\label{observer2}
\begin{cases}
\overset{.}{y}_0&=y_1\\
&~\vdots\\
\overset{.}{y}_n&=a_t-b_t\left(\sum_{i=0}^{n}k_i y_i-\sum_{i=2}^{n}k_i \epsilon^{n-i} \xi_i\right)\\
\overset{.}{\xi}&=\frac{1}{\epsilon}B \xi+ \begin{bmatrix}0\\\vdots\\0\\a_t-b_t\left(\sum_{i=0}^{n}k_i y_i-\sum_{i=2}^{n}k_i \epsilon^{n-i} \xi_i\right)\end{bmatrix}
\end{cases}
\end{equation}
where  $|a_t|\le L\left\|\overline{Y}(t)\right\|$ and $0<\underline{b}\le b_t \le \overline{b}$ for any $t\in[0,T)$, and where $B$ is a Hurwitz matrix. Then for any $(k_0,\cdots,k_n)\in \Omega_{L,\underline{b}, \overline{b},c}$, there exists $\epsilon^*>0$,  such that for any $0<\epsilon<\epsilon^*$,    $Y(t)$ and $\xi(t)$ satisfy
$$\|Y(t)\| \le ce^{-\beta t}(\|Y(0)\|+|y_0(0)|+\sqrt{2\lambda_{\max}(Q)}\|\xi(0)\|),$$
$$\|\xi(t)\|\le  \frac{c\left\|\overline{Y}(0)\right\|+\sqrt{\lambda_{\max}(Q)}\|\xi(0)\|}{\sqrt{\lambda_{\min}(Q)}}e^{-\beta t}, \forall t\in[0,T),$$ for some  $\beta>0$,
where $Q$ is the unique positive definite matrix satisfying $B^T Q+Q B=-I$.
\end{Lemma}

\begin{Remark}
Note that in Lemmas 2-4, we do not need the uncertain  functions $a_t$ and $b_t$ to have a fixed form, nor other conditions   excepts some upper bound functions, thanks to the design of the extended PID parameters and to the strong robustness of the PID controller. This enables us to establish global or semi-global results on output regulation of general SISO affine nonlinear uncertain systems by using its  (local) nonlinear canonical form in the analysis, as will be demonstrated shortly in the proofs of the theorems.
\end{Remark}

To prove our theorems, we also need the following result:

Theorem A1. Let $\Phi:\mathbb{R}^n\to \mathbb{R}^n$ be continuously differentiable. Then $\Phi$ is a global diffeomorphism if and only if the Jacobian matrix $J_\Phi(x)$ is nonsingular for all $x\in \mathbb{R}^n$ and $\lim_{\|x\|\to \infty}\|\Phi(x)\|=\infty$.

See \cite{sandberg1980global} and \cite{wu1972global} for detailed discussion.

\textbf{Proof of Theorem 1.}

Step 1. First, let us denote
\begin{equation}\label{222}
z(t)=(z_1(t),\cdots,z_n(t))^T\overset{\triangle}{=}\Phi(x(t))
\end{equation}
By Assumption (A), we know that the system (\ref{1}) has uniform relative degree $n$, thus   we have (see e.g., \cite{isidori2013nonlinear}),
\begin{equation}\label{3.00}
\begin{cases}
\overset{.}{z}_1&=z_2\\
&~\vdots\\
\overset{.}{z}_{n-1}&=z_n\\
\overset{.}{z}_n&=F(x)+G(x)u
\end{cases}
\end{equation}
where $z_1(t)=h(x(t))=y(t)$.
Denote $y_0(t)\overset{\triangle}{=}-\int_{0}^{t} e(s)ds-\frac{F(x^*)}{k_0 G(x^*)}$, $y_1(t)\overset{\triangle}{=}-e(t)=-(y^*-y(t))=z_1(t)-y^*$, $y_2(t)\overset{\triangle}{=}-\overset{.}{e}(t)=z_2(t)$, $\cdots$, $y_n(t)\overset{\triangle}{=}-e^{(n-1)}(t)=z_n(t)$. Then we have
\begin{equation}\label{closedddd}
\begin{cases}
\overset{.}{y}_0&=y_1\\
&~\vdots\\
\overset{.}{y}_{n-1}&=y_n\\
\overset{.}{y}_n&=F(x)+G(x)u\\
\end{cases}
\end{equation} where $y_0(0)=-\frac{F(x^*)}{k_0G(x^*)}$. Since $e(t)=y^*-y(t)=y^*-z_1(t)$,   we can rewrite  $u(t)$  as follows:
\begin{align}\label{upid} u(t)=&k_{1} e(t)+k_{0} \int_{0}^{t} e(s)ds+k_{2} \overset{.}{e}(t)+\cdots+k_{n} e^{(n-1)}(t)\nonumber\\=&-\bigg(\sum_{i=0}^{n}k_i y_i(t)+\frac{F(x^*)}{ G(x^*)}\bigg)\end{align}
Denote $a_t\overset{\bigtriangleup}{=} F(x(t))-\frac{F(x^*)}{ G(x^*)}G(x(t))$, $b_t\overset{\bigtriangleup}{=}G(x(t))$. From (\ref{closedddd})-(\ref{upid}), we have
\begin{equation}\label{closeddd}
\begin{cases}
\overset{.}{y}_0&=y_1\\
&~\vdots\\
\overset{.}{y}_{n-1}&=y_n\\
\overset{.}{y}_n&=a_t-b_t\sum_{i=0}^{n}k_i y_i.\\
\end{cases}
\end{equation}
Denote  $$Y(t)\overset{\triangle}{=}(y_1(t),\cdots, y_n(t))^T,~\overline{Y}(t)\overset{\triangle}{=}(y_0(t),\cdots,y_n(t))^T,$$  it is easy to see that $Y(t)=z(t)-z^*=\Phi(x(t))-\Phi(x^*)$, where $z^*$ is defined in (\ref{zzz}).

Step 2. Next, we will apply Lemma 3  to prove that if the initial state $\|x(0)\|\le R$ and the parameters $(k_0,\cdots,k_n)\in \Omega_{L_0,\underline{b}, b_0,c}$, where $b_0\overset{\triangle}{=}\tau_1(cR_0+|y^*|)$, $L_0\overset{\triangle}{=}\sup_{0\le\rho\le  R_0}\frac{\tau_1(|y^*|)+\underline{b}}{\underline{b}}\frac{\tau_2(\rho)}{\rho}$ and $R_0\overset{\triangle}{=}\tau_1(R)+|y^*|+\tau_1(|y^*|)$, then there exists $\alpha>0$, such that
 \begin{align}\label{yt}
 \left\|\overline{Y}(t)\right\| \le c e^{-\alpha t} \left\|\overline{Y}(0)\right\|,~~ t\in[0,T),
 \end{align} where $[0,T)$ is the maximal  existence interval of the  solution of  the closed-loop system (\ref{1}) and (\ref{6}) with initial state $x(0)$.

By Assumption (B), we have $\|H(x^*)\|\le \tau_1(\|\Phi(x^*)\|)=\tau_1(\|z^*\|)=\tau_1(|y^*|)$.  Therefore, $|F(x^*)|\le \|H(x^*)\|\le \tau_1(|y^*|)$ and we can obtain
\begin{align}\label{at}
&|a_t|=\left|F(x(t))-F(x^*) G(x(t))/ G(x^*)\right|\nonumber
\\\le&  |F(x(t))-F(x^*)|+|F(x^*)|  |G(x(t))-G(x^*)|/G(x^*)\nonumber
\\\le& |F(x(t))-F(x^*)|+\tau_1(|y^*|) |G(x(t))-G(x^*)|/\underline{b}\nonumber
\\\le& (\underline{b}+\tau_1(|y^*|)) \|H(x(t))-H(x^*)\|/\underline{b}
\le (\underline{b}+\tau_1(|y^*|)) \tau_2(\|\Phi(x(t))-\Phi(x^*)\|)/ \underline{b}\nonumber
\\=&(\underline{b}+\tau_1(|y^*|))\tau_2(\|Y(t)\|)/\underline{b}
\le(\underline{b}+\tau_1(|y^*|))\tau_2\left(\left \|\overline{Y}(t)\right \|\right)/\underline{b}.\end{align}

On the other hand, by Assumptions (A) and (B), it is easy to see that
\begin{align}\label{bt}
0<\underline{b}\le b_t=&G(x(t)) \le\|H(x(t))\|\le \tau_1(\|\Phi(x(t))\|)\nonumber \\ =&\tau_1(\|Y(t)+z^*\|)\le \tau_1(\|Y(t)\|+|y^*|)
\le \tau_1(\left \|\overline{Y}(t)\right \|+|y^*|).
\end{align}
Since $\|x(0)\|\le R$, then by Assumption (B), we have $$\|Y(0)\|=\|\Phi(x(0))-z^*\| \le \tau_1(\|x(0)\|)+\|z^*\|\le \tau_1(R)+|y^*|.$$   Recall that $y_0(0)=-\frac{F(x^*)}{k_0G(x^*)}$, we have
\begin{align*}\left\|\overline{Y}(0)\right\|\le& \|Y(0)\|+|y_0(0)|
\le  \tau_1(R)+|y^*|+\frac{\tau_1(|y^*|)}{k_0 G(x^*)}
\le \tau_1(R)+|y^*|+\tau_1(|y^*|),\end{align*}
where the last inequality holds since $k_0=\frac{\prod_{i=0}^{n}\lambda_{i}}{\underline{b}}\geq \frac{1}{\underline{b}}$ and $G(x^*)\geq \underline{b}$.
By Lemma 3, we know that (\ref{yt}) holds  for any  $(k_0,\cdots,k_n)\in \Omega_{L_0,\underline{b}, b_0,c}$ as long as the initial state $\|x(0)\|\le R$.

Step 3. In this step, we will show that if Assumption (C) holds,    then the maximal existence interval of the  solution of  the closed-loop system (\ref{1}) and (\ref{6}) is $[0,\infty)$.

 We use the contradiction argument.
Suppose that the solution of  the  system (\ref{1}) with controller (\ref{6}) only exists in $[0,T)$ for $T<\infty$ for some initial state $\|x(0)\|\le R$.
Denote  $x_0(t)=\int_{0}^{t}y^*-h(x(s))ds$, then $x_0(0)=0$  and by simple calculations, the extended PID controller (\ref{6}) can be rewritten as $u(t)=k_0 x_0(t)+k_1(y^*-h(x(t)))-k_2L_f h(x(t))-\cdots-k_n L_f^{n-1}h(x(t))$.    Therefore, the maximal existence interval of the solution of the following $(n+1)^{\text{th}}$ order autonomous  differential equation is also finite:
\begin{equation}{\label{e}}
\begin{cases}
\overset{.}{x}_0&=y^*-h(x)\\
\overset{.}{x}&=f(x)+g(x)(k_0x_0+k_1y^*-\sum_{j=1}^{n}k_jL_f^{j-1} h(x))
\end{cases}
\end{equation}
with the initial value $[0,x^T(0)]\in \mathbb{R}^{n+1}$. Then it is well-known from the theory of ordinary differential equations that \begin{align}\label{39}\limsup_{0\le t< T} \left\|[x_0(t),x(t)^T]\right\|=\infty.\end{align}
By step 2, we know that $ \left \|\overline{Y}(t)\right \|$ is bounded on $t\in[0,T)$. Hence, by (\ref{at}),(\ref{bt}), the right hand of (\ref{closeddd})  and the boundedness of $\left \|\overline{Y}(t)\right \|$, it is not difficult to conclude that
\begin{small}\begin{align*} \left\|\overset{.}{\overline{Y}}(t) \right\| \le N, ~~\forall t\in[0,T)\end{align*}\end{small}
 for some  constant $N>0$(possibly depends on the initial state $x(0)$).

 On the other hand, since $z(t)=\Phi(x(t))$, we have $\overset{.}{z}(t)=J_{\Phi}(x(t)) \overset{.}{x}(t)$ and therefore we can obtain
 \begin{align}\label{dot}
 \|\overset{.}{x}(t)\| &=\left\|J_\Phi ^{-1}(x(t)) \overset{.}{z}(t)\right\| \le \|J_\Phi ^{-1}(x(t))\|\| \overset{.}{z}(t)\|
 =\|J_\Phi ^{-1}(x(t))\| \left\|\overset{.}{Y}(t)\right\|   \nonumber\\
 &\le \| J_\Phi ^{-1}(x(t))\| \left\|\overset{.}{\overline{Y}}(t) \right\|
 \le N\left\|J_\Phi ^{-1}(x(t))\right\|,~~ t\in [0,T).
 \end{align}
By Assumption (C), we have \begin{align}\label{41}\|\overset{.}{x}(t)\| \le \alpha_1\|x(t)\|\log \|x(t)\|+\alpha_2 \end{align} for any $t\in [0,T)$, where $\alpha_1=NN_1$ and $\alpha_2=NN_2$. Denote $v(t)\overset{\triangle}{=}\|x(t)\|$ and $D^+ v(t)\overset{\triangle}{=}\limsup_{h\to 0^+} \frac{v(t+h)-v(t)}{h}$ be the upper right-hand derivative of $v(t).$

Then it is not difficult to obtain \begin{align*}
D^+ v(t)\le &|D^+ v(t)|=\limsup_{h\to 0^+} \big|\frac{\|x(t+h)\|-\|x(t)\|}{h}\big|
\le  \limsup_{h\to 0^+} \frac{\|x(t+h)-x(t)\|}{h}\\=& \limsup_{h\to 0^+} \big\|\frac{x(t+h)-x(t)}{h}\big\|= \|\overset{.}{x}(t)\|\end{align*}
Noticing that $v(t)=\|x(t)\|$, from (\ref{41}), we have
\begin{align*}
D^+ v(t) \le \|\overset{.}{x}(t)\|\le \alpha_1v(t)\log v(t)+\alpha_2,~~t\in[0,T)
\end{align*}
By the comparison lemma in ordinary differential equations(see e.g., \cite{khalil2002nonlinear}),  we have $$\int_{v(0)}^{v(t)}\frac{d\eta}{\alpha_1 \eta \log\eta +\alpha_2}\le t< T, ~~\forall t\in [0,T),$$ which implies \begin{align}\label{xt}\sup_{0\le t< T} v(t)=\sup_{0\le t< T} \|x(t)\|<\infty\end{align} since $\int_{v(0)}^{\infty}\frac{d\eta}{\alpha_1 \eta \log\eta +\alpha_2}=\infty$. By (\ref{xt}) and the fact $T<\infty$, it is not difficult to see that $$\sup_{0\le t< T} |x_0(t)|= \sup_{0\le t< T}  \left|\int_{0}^{t}y^*-h(x(s))ds\right|<\infty.$$ Therefore, the solution of (\ref{e}) with initial state  $[0,x^T(0)]$ satisfy $$\sup_{0\le t< T} \left\|[x_0(t),x^T(t)]\right\|<\infty,$$  which contradicts to (\ref{39}).

Therefore, under Assumptions (A),(B) and (C), if $(k_0,\cdots,k_n)\in \Omega_{L_0,\underline{b}, b_0,c}$, then for any initial state $\|x(0)\|\le R$, the solution of the closed-loop system will exist in $[0,\infty)$.

Step 4. If Assumption (C) is replaced by Assumption (C'),    then the maximal existence interval of the  solution of  the closed-loop system (\ref{1}) and (\ref{6}) is also $[0,\infty)$.

 We use the contradiction argument again.
Suppose that the solution of  the  system (\ref{1}) with controller (\ref{6}) only exists in $[0,T)$ for $T<\infty$ and for some initial state $\|x(0)\|\le R$.
By step 2, we know that $Y(t)$ is bounded on $[0,T)$. Notice that $Y(t)=\Phi(x(t))-z^*$, we know that $\Phi(x(t))$ is bounded on $[0,T)$. Since $\Phi$ is a global diffemorphism on $\mathbb{R}^n$, by Theorem (A1), we know that $x(t)$ is bounded on $[0,T)$. Similarly, we have $\sup_{0\le t< T} \left\|[x_0(t),x(t)^T]\right\|<\infty,$ which contradicts to (\ref{39}).

Step 5. Since  solution of the closed-loop equation exists in $[0,\infty)$, we conclude that (\ref{at}) and (\ref{bt}) are satisfied in $[0,\infty)$. By using Lemma 3 again, we have $\left\|\overline{Y}(t)\right\| \le c e^{-\alpha t} \left\|\overline{Y}(0)\right\|$ for any $t\in[0,,\infty)$.

Therefore, we have
\begin{small}
\begin{align*}|e(t)|&=|y_1(t)| \le \|Y(t)\|=\|\Phi(x(t))-z^*\|  \le \left\|\overline{Y}(t)\right\| \\
&\le c e^{-\alpha t} \left\|\overline{Y}(0)\right\|\le  c e^{-\alpha t} \left(\|Y(0)\|+\left|\frac{F(x^*)}{k_0G(x^*)}\right|\right)\\
&\le c e^{-\alpha t} \left(\|\Phi(x(0))-z^*\|+\left|\frac{F(x^*)}{k_0G(x^*)}\right|\right)\\
&= c e^{-\alpha t} \left(\|\Phi(x(0))-z^*\|+\frac{\tau_1(|y^*|)}{k_0\underline{b}}\right)\end{align*}\end{small} for any $t\in [0,\infty).$

If Assumption (C') is satisfied, then we can see that $\Phi^{-1}(z^*)$ only has one element, denote it as $x^*$. It is not difficult to obtain $\lim_{t\to \infty}x(t)= x^*$. As a consequence,  the  system state $x(t)$ is bounded.
This completes the proof of Theorem 1. $\Box$

\textbf{Proof of Corollary 1.}

By Assumption (A), we know that $J_\Phi(x)$ is invertible for any $x\in\mathbb{R}^n$. From Assumption (B0), we have $\lim_{\|x\|\to \infty}\|\Phi(x)\|\geq \lim_{\|x\|\to \infty}\rho_0(\|x\|)=\infty$.  By using Theorem A1, we know that $\Phi$ is a global diffeomorphism on $\mathbb{R}^n$, i.e., Assumption (C') is satisfied.  Let $y^*$ be any  given setpoint, then  $\Phi^{-1}(z^*)$ is not empty  and $\Phi^{-1}(z^*)$  only has one element, we denote it as $x^*$.

Now, we proceed to verify Assumption (B), i.e., to find two increasing  functions $\tau_1,\tau_2$ with  $\limsup_{\rho\to0}  \frac{\tau_2(\rho)}{\rho}< \infty$ such that the following inequalities  \begin{align}\label{MH}
\|\Phi(x)\|\le \tau_1(\|x\|),~\|H(x)\|\le \tau_1(\|\Phi(x)\|),\end{align}
\begin{align}\label{HM}\|H(x)-H(x^*) \| \le \tau_2(\|\Phi(x)-\Phi(x^*)\|)\end{align} hold for any $x\in\mathbb{R}^n$.

Since $\rho_0(\|x\|)\le \|\Phi(x)\|$, it is easy to obtain \begin{align}\label{xx}\|x\|\le \rho_0^{-1}(\|\Phi(x)\|), \forall x\in\mathbb{R}^n,\end{align} where $\rho_0^{-1}(r)\overset{\triangle}{=}\sup\{y \geq 0 | \rho_0(y)\le r\}$ for $r\geq 0.$  From (\ref{xx}), we have
  \begin{align}\label{zz}\|\Phi^{-1}(z)\|\le \rho_0^{-1}(\|z\|), \forall z\in\mathbb{R}^n.\end{align}
By Assumption (B0) and (\ref{xx}), we have
\begin{align}\label{H}\|H(x)\|\le \rho_1(\|x\|)\le  \rho_1\circ \rho_0^{-1}(\|\Phi(x)\|),\end{align} where $ \rho_1\circ \rho_0^{-1}$ denotes the composition of the functions $\rho_1$ and $ \rho_0^{-1}$. Therefore, we conclude that (\ref{MH}) holds with $$\tau_1(r)=\max\{\rho_1(r),\rho_1\circ \rho_0^{-1}(r)\}.$$

From (\ref{xx}),  we have  $\|x^*\|\le \rho_0^{-1}(\|\Phi(x^*)\|)=\rho_0^{-1}(\|z^*\|)\le \rho_0^{-1}(|y^*|)\overset{\triangle}{=}M_0$.
Next, we proceed to estimate the upper bound of $\|H(x)-H(x^*)\|$.
For this, we need to prove the following statement first: Let $U\subset \mathbb{R}^n$ be a convex open set. $f:U\to \mathbb{R}^m$ is a continuously differentiable vector valued function. Then for any $x_1,x_2\in U,$ we have
  \begin{align}\label{46}
\|f(x_1)-f(x_2)\|\le \sup_{0\le \theta \le 1}\left\|J_f( x_2+\theta(x_1-x_2)\right\|\|x_1-x_2\|\end{align}
The proof is elementary. Since $U$ is convex, then $g(\theta)=f(x_2+\theta(x_1-x_2))$ is a continuously differentiable  function of $0\le \theta\le 1.$ Note that
$f(x_1)-f(x_2)=g(1)-g(0)=\int_{0}^{1}g'(\theta)d\theta=\int_{0}^{1}J_f(x_2+\theta(x_1-x_2))(x_1-x_2)d\theta .$
Thus, we have
\begin{align*}
\|f(x_1)-f(x_2)\|\le& \int_{0}^{1}\left\|J_f(x_2+\theta(x_1-x_2))\right\|d\theta \|x_1-x_2\|\\
\le &\sup_{0\le \theta \le 1}\left\|J_f(x_2+\theta(x_1-x_2))\right\|\|x_1-x_2\|
\end{align*}

By (\ref{46}), we know that
$$\|H(x)-H(x^*)\|\le\sup_{0\le \theta \le1} \left\|J_H(x^*+\theta(x-x^*))\right\|\|x-x^*\|.$$ Since $\left\|J_H(x)\right\|\le \rho_1(\|x\|)$, $\|x^*\|\le M_0$ and $0\le \theta\le 1$,  therefore we have
\begin{align}\label{H1} \|H(x)-H(x^*)\| \le \rho_1(M_0+\|x-x^*\|)\|x-x^*\|.\end{align}
On the other hand, from the identity $\Phi^{-1}\circ\Phi(x)=x$, we have $J_{\Phi^{-1}}(\Phi(x))J_ \Phi(x)=I$,  i.e., \begin{align}\label{phi}
J_{\Phi^{-1}}(\Phi(x))=J_ \Phi^{-1}(x),~~ \forall x\in\mathbb{R}^n,
\end{align}where $J_{\Phi^{-1}}$ denotes the Jacobian matrix of the mapping $\Phi^{-1}$.
By using (\ref{46}) again, we have
 \begin{align}\label{33}
 \|x-x^*\|=&\|\Phi^{-1}(\Phi(x))-\Phi^{-1}(\Phi(x^*))\|\nonumber\\
 \le& \sup_{0\le \theta \le 1}\left\|J_{\Phi^{-1}} (z_\theta) \right\|\|\Phi(x)-\Phi(x^*)\|\end{align}
 where $z_\theta \overset{\triangle}{=}\theta (\Phi(x)-\Phi(x^*))+\Phi(x^*)$.
 Combine  (\ref{phi}) and (\ref{33}), we have
  \begin{align}\label{x,x}
  \|x-x^*\| \le \sup_{0\le \theta \le 1}\left\|J_ \Phi^{-1}(\Phi^{-1}(z_\theta)) \right\|\|\Phi(x)-\Phi(x^*)\|\end{align}  Since $\|z_\theta\|\le \|\Phi(x)-\Phi(x^*)\|+\|\Phi(x^*)\| =\|\Phi(x)-\Phi(x^*)\|+|y^*|$, from (\ref{zz}), we have
  \begin{align}\label{36}\|\Phi^{-1}(z_\theta)\|\le \rho_0^{-1}(\|\Phi(x)-\Phi(x^*)\|+|y^*|).\end{align}
By Assumption (B0), we have $\left\|J_\Phi ^{-1}(x)\right\| \le \rho_1(\|x\|)$, combine (\ref{x,x}) and (\ref{36}), we can obtain
\begin{align}
\label{x2}\|x-x^*\|\le \rho_2(\|\Phi(x)-\Phi(x^*)\|),
\end{align}
where $\rho_2(r)\overset{\triangle}{=}\left[\rho_1\circ \rho_0^{-1}(r+|y^*|) \right] r .$

Denote $\tau_2(r) \overset{\triangle}{=} \rho_1(M_0+\rho_2(r))\rho_2(r)$, then from (\ref{H1}) and (\ref{x2}), we have
 \begin{align} \|H(x)-H(x^*)\| \le \tau_2(\|\Phi(x)-\Phi(x^*)\|).\end{align}
Since $$\limsup_{r \to 0}\frac{ \rho_2(r)}{r}=\limsup_{r \to 0} \rho_1\circ \rho_0^{-1}(r+|y^*|)<\infty,$$ we arrive at $\limsup_{r \to 0}\frac{ \tau_2(r)}{r}<\infty.$

Therefore, all the conditions in Assumption (B) are satisfied. Consequently, all results in  Theorem 1 hold.

Finally, we show that $x(t)$ has exponential convergence rate.
By the proof of Theorem 1, we  know that $\|\Phi(x(t))-z^*\| \le c e^{-\alpha t} \big(\|\Phi(x(0))-z^*\|+\frac{\tau_1(|y^*|)}{k_0\underline{b}}\big).$
As a consequence, $\Phi(x(t))$ will converge to $z^*$ exponentially.
From (\ref{x2}), we have $\|x(t)-x^*\|\le \rho_2(\|\Phi(x(t))-z^*\|)$. Recall that $\limsup_{r \to 0}\frac{ \rho_2(r)}{r}<\infty,$ we see that the trajectory $x(t)$ is bounded and $\lim_{t\to\infty} x(t)=x^*$ exponentially.$\Box$

\textbf{Proof of Theorem 2.} We will use the same notations as in the proof of Theorem 1, then
\begin{equation}\label{closedd}
\begin{cases}
\overset{.}{y}_0&=y_1\\
&~\vdots\\
\overset{.}{y}_n&=a_t-b_t(k_0y_0+\cdots+k_n y_n),
\end{cases}
\end{equation} where $a_t=F(x(t))-\frac{F(x^*)}{G(x^*)} G(x(t))$, $b_t=G(x(t))$ and $y_0(0)=-\frac{F(x^*)}{k_0G(x^*)}$.

By Assumption (B1),  we have $|F(x^*)|\le \tau(\|\Phi(x^*)\|)=\tau(|y^*|)$. Denote $\tilde{L}=\frac{L(\tau(|y^*|)+\underline{b})}{\underline{b}}$. By  $\|H(x)-H(x^*)\|\le L\|\Phi(x)-\Phi(x^*)\|$, we can obtain
\begin{align}\label{atl}|a_t|\le& |F(x(t))-F(x^*)|+\left|\frac{F(x^*)}{G(x^*)}\right| |G(x(t)-G(x^*)|\nonumber\\
\le & \frac{\tau(|y^*|)+\underline{b}}{\underline{b}} \|H(x(t))-H(x^*)\|\le \tilde{L}\|\Phi(x(t))-z^*\|\nonumber\\
=&\tilde{L}\|z(t)-z^*\|=\tilde{L}\|Y(t)\|\le \tilde{L}\left \|\overline{Y}(t)\right \|,\end{align}  and from $0<\underline{b}\le G(x)\le \overline{b}$, we have $$0<\underline{b}\le b_t=G(x(t))\le \overline{b}.$$
By Lemma 2, we know that $ \left \|\overline{Y}(t)\right \| \le  c  \left \|\overline{Y}(0)\right \|$ for $t\in[0,T)$ if the parameters $(k_0,\cdots,k_n)\in \Omega_{\tilde{L},\underline{b}, \overline{b},c}$, where $[0,T)$ is the maximal  existence interval of the  solution of  the closed-loop system (\ref{1}) and (\ref{6}).

Similar to the proof of Theorem 1,  we know  that either Assumption (C) or (C') holds, then  the solution of the closed-loop system will exist in $[0,\infty)$ for any initial state $x(0)\in \mathbb{R}^n$, whenever $(k_0,\cdots,k_n)\in \Omega_{\tilde{L},\underline{b}, \overline{b},c}$. By using Lemma 2 with $T=\infty$, we have $\left\|\overline{Y}(t)\right\| \le c e^{-\alpha t} \left\|\overline{Y}(0)\right\|$ for any $t\geq 0$ and for some  $\alpha>0$. Therefore, $$|e(t)|\le \left \|\overline{Y}(t)\right\| \le c e^{-\alpha t} \|\overline{Y}(0)\|\le c e^{-\alpha t} (\|Y(0)\|+\frac{\tau(|y^*|)}{k_0\underline{b}}),$$ for any $t\in [0,\infty)$, which implies $e(t)$ converge to $0$ exponentially.

The proof of the convergence of the state $x(t)$ under Assumption (C') is similar to that of Theorem 1, which will not be repeated .$\Box$

\textbf{Proof of Theorem 3.}

Denote $y_0(t)=-\int_{0}^{t} e(s)ds-\frac{F(x^*)}{k_0G(x^*)}$, $y_1(t)=-e(t)=z_1(t)-y^*$, $\cdots$, $y_n(t)=-e^{(n-1)}(t)=z_n(t)$, $Y(t)=(y_1(t),\cdots,y_n(t))^T=z(t)-z^*$, $\xi_i(t)=\frac{y_i(t)+\hat{z}_i(t)}{\epsilon^{n-i}}$, $i=1,\cdots,n$. Then the differential observer-based extended PID controller  (\ref{5})-(\ref{4}) can be rewritten as
\begin{align*}u(t)=&k_{1} e(t)+k_{0} \int_{0}^{t} e(s)ds+k_{2} \hat{z}_2(t)+\cdots+k_{n} \hat{z}_n(t)\\=&-\sum_{i=0}^{n}k_i y_i(t)+\sum_{i=2}^{n}k_i \epsilon^{n-i} \xi_i(t)-\frac{F(x^*)}{G(x^*)}.\end{align*}
 Consequently, we  have the following equation:
\begin{equation}\label{observer}
\begin{cases}
\overset{.}{y}_0&=y_1\\
&~\vdots\\
\overset{.}{y}_{n-1}&=y_n\\
\overset{.}{y}_n&=a_t-b_t\left(\sum_{i=0}^{n}k_i y_i-\sum_{i=2}^{n}k_i \epsilon^{n-i} \xi_i\right)\\
\overset{.}{\xi}_1&=-\frac{\beta_1}{\epsilon}\xi_1+\frac{\xi_2}{\epsilon}\\
&~\vdots\\
\overset{.}{\xi}_{n-1}&=-\frac{\beta_{n-1}}{\epsilon}\xi_1+\frac{\xi_n}{\epsilon}\\
\overset{.}{\xi}_{n}&=-\frac{\beta_{n}}{\epsilon} \xi_1+a_t-b_t\left(\sum_{i=0}^{n}k_i y_i-\sum_{i=2}^{n}k_i \epsilon^{n-i} \xi_i\right)
\end{cases}
\end{equation} where $a_t=F(x(t))-\frac{F(x^*)}{G(x^*)} G(x(t))$, $b_t=G(x(t))$.

Denote $B\overset{\triangle}{=}\begin{bmatrix}-\beta_1&1&0&\cdots&0\\\vdots&&\ddots&&\vdots\\-\beta_{n-1}&0&0&\cdots&1\\-\beta_n&0&0&\cdots&0\\
\end{bmatrix},$ and $\xi(t)\overset{\triangle}{=}(\xi_1(t),\cdots,\xi_n(t))^T$. Then (\ref{observer}) turns into a more compact form:
\begin{equation}\label{observer1}
\begin{cases}
\overset{.}{y}_0=y_1\\
~~~~\vdots\\
\overset{.}{y}_n=a_t-b_t\left(\sum_{i=0}^{n}k_i y_i-\sum_{i=2}^{n}k_i \epsilon^{n-i} \xi_i\right)\\
\overset{.}{\xi}=\frac{1}{\epsilon}B \xi+ \left(0,\cdots,a_t-b_t\left(\sum_{i=0}^{n}k_i y_i-\sum_{i=2}^{n}k_i \epsilon^{n-i} \xi_i\right)\right)^T
\end{cases}
\end{equation}
By Assumption (B1), and from the proof of Theorem 2, we have
$$|a_t|\le \tilde{L}\|Y(t)\|,~\underline{b}\le b_t=G(x(t))\le \overline{b},$$ where $\tilde{L}=\frac{L(\tau(|y^*|)+\underline{b})}{\underline{b}}$. From Lemma 4, we know that for any $(k_0,\cdots,k_n)\in \Omega_{\tilde{L},\underline{b}, \overline{b},c}$, there exists $\epsilon^*>0$, such that for any  $0<\epsilon< \epsilon^*$,   there exists $\beta>0$, such that  the following inequality holds on the interval where the solution exists:
$$\|Y(t)\| \le ce^{-\beta t}(\|Y(0)\|+|y_0(0)|+\sqrt{2\lambda_{\max}(Q)}\|\xi(0)\|).$$
The final thing we need to prove is: if the controller parameters $(k_0,\cdots,k_n)\in \Omega_{\tilde{L},\underline{b}, \overline{b},c}$ and the observer gain $0<\epsilon< \epsilon^*$,   then for any initial states $x(0)\in\mathbb{R}^n$, $\hat{z}(0)\in \mathbb{R}^n$, the maximal existence interval of the  solution of  the closed-loop system (\ref{1}),(\ref{5}) and (\ref{4}) is $[0,\infty)$.

We use the contradiction argument. Suppose that for some initial states $x(0)\in\mathbb{R}^n$, $\hat{z}(0)\in \mathbb{R}^n$, the  solution of  the closed-loop system  only exists in $[0,T)$ for some $T<\infty$.
Then we can  see  that the maximal interval of existence   of the following $(2n+1)^{\text{th}}$-order autonomous  differential equation composed of the control systems and the observers is also finite:
\begin{equation}\label{de}
\begin{cases}
\overset{.}{x}_0&=y^*-h(x)\\
\overset{.}{x}&=f(x)+g(x)\big(k_{0} x_0+k_{1} (y^*-h(x))+\sum_{j=2}^nk_{j} \hat{z}_j\big)\\
\overset{.}{\hat{z}}_1&=\hat{z}_2+\frac{\beta_1}{\epsilon}(y^*-h(x)-\hat{z}_1)\\
&\vdots\\
\overset{.}{\hat{z}}_n&=\frac{\beta_n}{\epsilon^n}(y^*-h(x)-\hat{z}_1)
\end{cases}
\end{equation}
for   the initial conditions $\begin{bmatrix}0&x(0)^T&\hat{z}(0)^T\end{bmatrix}\in \mathbb{R}^{2n+1}$. Therefore, the solution of (\ref{de}) with the initial value $\begin{bmatrix}0&x(0)^T&\hat{z}(0)^T\end{bmatrix}$ will satisfy
\begin{align}\label{41.0}\limsup_{t\to T}\|\begin{bmatrix}x_0(t)&x(t)^T&\hat{z}(t)^T\end{bmatrix}\|=\infty.\end{align}
From Lemma 4, we know that $\left\|\overline{Y}(t)\right\|$ and $\|\xi(t)\| $ are bounded on $[0,T)$. Therefore, by the right hand  of (\ref{observer1}), we can obtain $$\left\|\overset{.}{Y}(t)\right\|\le N, ~~\forall t\in[0,T)$$
for some $N>0$.  Similar to the proof of Theorem 1,  we know that   $\limsup_{t\to T}\|x(t)\| <\infty$. As a consequence, $x_0(t)=\int_{0}^{t}y^*-h(x(s))ds$ is also bounded on $[0,T)$. On the other hand, by the boundedness of $\xi(t)$, we have $\sup_{0\le t< T} \|\hat{z}(t)\|<\infty$,  which  contradicts to (\ref{41.0}).

Finally, by noticing $|y_0(0)|\le \frac{\tau(|y^*|)}{k_0 \underline{b}}$ and applying Lemma 4 again, we conclude that \begin{align*}&|e(t)|\le \|Y(t)\|\le ce^{-\beta t}\left(\|\Phi(x(0))-z^*\|+\frac{\tau(|y^*|)}{k_0\underline{b}}+\sqrt{2\lambda_{\max}(Q)}\|\xi(0)\|\right)\end{align*}

The proof of the convergence of the state $x(t)$ under Assumption (C') is similar to that of Theorem 1, which will not be repeated .$\Box$

\section{Conclusions}
In this paper, we have presented a theoretical investigation on the extended PID controller for a general class of SISO affine-nonlinear uncertain dynamical systems.  It is shown that the extended PID controller  can globally (or semi-globally) stabilize this class of systems  under some fairly general   conditions on  nonlinearity  and uncertainty of the systems, as long as the controller parameters are chosen from an open and unbounded  parameter manifold. It is worth noting that   the extended PID controller includes the classical PID controller as a special case, and that its  design  does not need the precise information  about the mathematical model. It is also interesting to note that  the  nonlinear canonical form (which may not be global) in  geometrical nonlinear control theory can be used in our theoretical analysis to get global(or semi-global) convergence results, thanks to the strong robustness of the extended PID controller. This enables us to avoid assuming special system structures like pure-feedback forms and to get general conditions that can be considerably simplified once the system structures are in certain special forms.   Of course,  many interesting problems still remain open. It would be interesting to consider extended  PID control for  multi-input-multi-output affine nonlinear uncertain systems, and  to generalize our recent results on PID control of coupled multi-agent dynamical systems   \cite{yuan2018uncoupled}.    It would also be interesting to consider more complicated situations such as saturation, deadzone, time-delayed inputs,  sampled-data PID controllers under a prescribed sampling rate, etc. These belong to further investigation.

\section{Appendix}

\subsection{Appendix A}
\textbf{Proof of Lemma 1.} From the definition of $\Omega_1$,  it is easy to see that $\lambda_n>2n+2>2n+1>\lambda_{n-1}>\cdots>\lambda_0>2$ and $|\lambda_j-\lambda_i|\geq 1$ for $i\neq j$ when $(\lambda_0,\cdots,\lambda_{n})\in \Omega_1$.

First, we  prove that  $P$ is invertible and $(\Pi_{i=0}^n\lambda_i)^n\det{P}=\Pi_{0\le i<j\le n}(\lambda_i-\lambda_j)$.

By the properties of determinant, we have $(\Pi_{i=0}^n\lambda_i)^n\det{P}=\det {P'}$, where $P'$ is a Vandermonde matrix defined by
 \begin{align*}\label{P'}
  P'=\begin{bmatrix}1&\cdots&1\\-\lambda_0&\cdots&-\lambda_{n}\\
\vdots&&\vdots\\\\(-\lambda_{0})^{n}&\cdots&(-\lambda_{n})^{-n}
\end{bmatrix}
\end{align*}
Since $P'$ is a Vandermonde matrix, we have $(\Pi_{i=0}^n\lambda_i)^n\det{P}=\Pi_{0\le i<j\le n}(\lambda_i-\lambda_j)\neq 0$, which implies that $P$ is invertible.

 Next, we proceed to prove that there exists $\delta>0$, such that $|\det{P}|\geq \delta$ when $(\lambda_0,\cdots,\lambda_{n})\in \Omega_1$.

It is easy to see that
\begin{small}
\begin{align*}
\lambda_n^n (\Pi_{i=0}^{n-1}\lambda_i)^n|\det{P}|=&\Pi_{0\le i<j\le n-1}(\lambda_j-\lambda_i)\Pi_{i=0}^{n-1} (\lambda_n-\lambda_i)\\\geq&\Pi_{i=0}^{n-1} (\lambda_n-\lambda_i)\geq \Pi_{i=0}^{n-1} (\lambda_n-(2i+3))\\\geq &\lambda_n^n \Pi_{i=0}^{n-1} (1-\frac{2i+3}{\lambda_n})\geq \lambda_n^n (\frac{1}{2n+2})^n
\end{align*}
\end{small}
 Since  $(\Pi_{i=0}^{n-1}\lambda_i)^n$ is bounded, we see that  there exists $ \delta>0$ depends on $n$ only, such that $|\det{P}|\geq \delta$.

 By using the fact $|\det{P}|\geq \delta$ and the formula $P^{-1}=\frac{P^*}{\det{P}}$, where $P^*$ is the adjoint matrix of $P$, we conclude that $\|P^{-1}\|\le \frac{\|P^*\|}{\delta}$.

Since all $(n+1)\times (n+1)$ elements of $P$ are bounded (bounded by $1$),  it is easy to see that there exists $\delta_2$, such that   $\|P\| \le \delta_2$ and $\|P^*\|\le \delta_2$ for any $\lambda\in\Omega_1$.

Hence, $c_1=\sup_{\lambda\in \Omega_1}\|P\|<\infty$, and $c_2=\sup_{\lambda\in \Omega_1}\|P\|\|P^{-1}\|\le c_1 \frac{\delta_2}{\delta}<\infty.$

Recall that $(d_0,\cdots,d_n)^T$ is the last column of $P^{-1}$, by some simple calculations of determinants, we have $d_i=\frac{\lambda_i^n}{\prod_{j\neq i}(\lambda_i-\lambda_j)}$, $i=0,\cdots,n$.
Since $\lambda\in \Omega_1$, we have $$0<d_n=\frac{\lambda_n^n}{\prod_{i\neq n}(\lambda_n-\lambda_i)}\le \big(\frac{\lambda_n}{\lambda_n-(2n+1)}\big)^n<(2n+2)^n,$$ therefore  $c_3<\infty.$
For $i=0,\cdots,n-1$, it is easy to see that $$\frac{\lambda_n\lambda_i^n }{(\lambda_n-\lambda_i) \left|\prod_{j\neq i, j\neq n}(\lambda_j-\lambda_i)\right|}\le \frac{(2n+1)^n \lambda_n}{\lambda_n-(2n+1)},$$ and so $c_4(i)<\infty.$
Hence, $c_0<\infty$ and therefore the proof of Lemma 1 is complete.

\textbf{Proof of Lemma 2.}
 Rewrite (\ref{close}) as \begin{equation}\label{close1}
\begin{cases}
\overset{.}{y}_0=y_1\\
~~~~\vdots\\
\overset{.}{y}_n=-\underline{b}\sum_{i=0}^n k_iy_i+a_t+(\underline{b}-b_t)\sum_{i=0}^n k_iy_i
\end{cases}.
\end{equation}
Suppose that $(k_0,\cdots,k_n)\in \Omega_{L,\underline{b}, \overline{b},c}$ and denote $$A\overset{\triangle}{=}\begin{bmatrix}0&1&0&\cdots&0\\ \vdots&\vdots&\vdots&\ddots&\vdots\\0&0&0&\cdots&1\\-\underline{b}k_0&-\underline{b}k_1&-\underline{b}k_2&\cdots&-\underline{b}k_n\end{bmatrix}.$$
Then $(\ref{close1})$ can be rewritten as \begin{equation}\label{close11}
\overset{.}{\overline{Y}}=A \overline{Y} +(0,\cdots,0, a_t+(\underline{b}-b_t)(k_0y_0+\cdots+k_ny_n))^T
\end{equation}

It is easy to see that the characteristic polynomial of $A$ is $\det(sI-A)=s^{n+1}+\sum_{i=0}^{n} \underline{b} k_i s^i$.

By the definition (\ref{Omega}) of $\Omega_{L,\underline{b}, \overline{b},c}$, there exists $(\lambda_0,\cdots,\lambda_n)\in \Omega_{\Lambda }$ such that $( \underline{b} k_0, \cdots,   \underline{b} k_n)=(\prod_{i=0}^{n}\lambda_{i},\cdots,\sum_{i=0}^{n}\lambda_{i})$.
 Therefore, by Vieta's formulas, we know that $-\lambda_i, i=0,\cdots,n$ are $(n+1)$ distinct eigenvalues of $A$.
Hence,  $A$ is similar to $J$, where $J$ is a diagonal matrix defined by $J\overset{\triangle}{=}\textrm{diag}(-\lambda_0
,\cdots,-\lambda_{n})$.

It is not difficult to get the relationship $A P=P J$, where $P$ is defined in (\ref{P}). To simplify the analysis, we introduce an invertible linear transformation $\overline{Y}(t)=P\overline{w}(t)$, where $\overline{w}=(w_0,\cdots,w_{n-1},w_n)^T$.

We first need to prove the following equality:  $$  \sum_{i=0}^n k_ny_n=\frac{1}{\underline{b}}\sum_{i=0}^{n}\lambda_i w_i.$$
 By the relationship $\overline{Y}=P\overline{w}$, we have $y_i=\sum_{k=0}^{n} (-\lambda_k)^{i-n}w_k$.
Therefore
\begin{align*} \sum_{i=0}^{n}\underline{b} k_iy_i =&\sum_{i=0}^{n} \underline{b} k_i \sum_{k=0}^{n} (-\lambda_k)^{i-n}w_k
=\sum_{k=0}^{n}  \big(\sum_{i=0}^{n} \underline{b} k_i(-\lambda_k)^{i-n}\big)w_k \\
=&\sum_{k=0}^{n}  \big(\frac{\sum_{i=0}^{n} \underline{b} k_i(-\lambda_k)^{i}}{(-\lambda_k)^{n}}\big)w_k
=\sum_{k=0}^{n}  \frac{-(-\lambda_k)^{n+1}}{(-\lambda_k)^{n}}w_k=\sum_{k=0}^{n}\lambda_k w_k.
\end{align*} The second to last equality holds since $-\lambda_k$ is the root of the polynomial
$s^{n+1}+\sum_{i=0}^{n} \underline{b} k_i s^i$.

By the relationship   $A=P J P^{-1}$,
  it is easy to see $(\ref{close11})$ can be transformed to \begin{equation}\label{close1111}
\overset{.}{\overline{w}}=J \overline{w} +P^{-1}\bigg(0,\cdots,0, a_t+(\underline{b}-b_t)\sum_{i=0}^nk_iy_i\bigg)^T
\end{equation}
Since $(d_0,\cdots,d_n)^T$ is the last column of the matrix $P^{-1}$,  by  the equality $\sum_{i=0}^nk_iy_i=\frac{1}{\underline{b}} \sum_{i=0}^{n}\lambda_i w_i$, we see that  (\ref{close1111}) becomes
\begin{equation}\label{8}
\begin{cases}
\overset{.}{w}_0=-\lambda_0 w_0+d_0 \left( a_t+\frac{\underline{b}-b_t}{\underline{b}}
\sum_{i=0}^{n}\lambda_i w_i\right) \\
~~~~~\vdots\\
\overset{.}{w}_{n}=-\lambda_{n} w_{n}+d_{n} \left (a_t+\frac{\underline{b}-b_t}{\underline{b}}
\sum_{i=0}^{n}\lambda_i w_i\right).
\end{cases}
\end{equation}

Now, we consider the following quadratic function: $$V(\overline{w}(t))=\frac{1}{2}\sum_{i=0}^{n}w_i^2(t)=\frac{1}{2}\|\overline{w}(t)\|^2.$$

Then it is easy to compute the time derivative of $V$  $\overset{.}{V}(\overline{w}(t))\overset{\triangle}{=}\frac{d V(\overline{w}(t))}{dt}$ as follows:
\begin{align}{\label{V}}
\overset{.}{V}(\overline{w}(t))=&-\sum_{i=0}^{n}\lambda_i w_i^2+\bigg(\sum_{i=0}^{n}d_i w_i\bigg)\bigg(a_t+\frac{\underline{b}-b_t}{\underline{b}}\sum_{i=0}^{n}\lambda_i w_i\bigg)\nonumber\\
=&\underbrace{-\sum_{i=0}^{n}\lambda_i w_i^2}_{\mathrm{I}}+\underbrace{\bigg(\sum_{i=0}^{n}d_i w_i\bigg)a_t}_{\mathrm{II}}
+\underbrace{\frac{\underline{b}-b_t}{\underline{b}}\sum_{i=0}^{n} d_i w_i \sum_{i=0}^{n}\lambda_i w_i}_{\mathrm{III}}
\end{align}
Next, we proceed to estimate  (\ref{V}) term by term.

Denote $(w_0,\cdots,w_{n-1})^T\overset{\triangle}{=}w$. Obviously, the first term \begin{equation}\label{I}\mathrm{I}= -\sum_{i=0}^{n}\lambda_i w_i^2 \le -2\|w\|^2-\lambda_n w_n^2 \end{equation} since $\lambda_i> 2,i=0\cdots,n-1$.

By Lemma 1, we have $\|P\|\le c$, therefore $|a_t|\le L \left \|\overline{Y}\right\|=  L \left\|P \overline{w}\right\| \le L c \|\overline{w}\|\le L c(\|w\|+|w_n|)$.

On the other hand,  by Lemma 1 and the fact $c\geq c_0$, we also have $|d_i|\le \frac{c}{(2n+1) n\lambda_{n}}<  \frac{c}{\sqrt{n}\lambda_{n}}, i=0,\cdots, n-1 $, and $|d_n|\le \frac{c}{(2n+1)\sqrt{n}}<c$. Therefore, we have $|\sum_{i=0}^{n}d_i w_i| \le c \big(\sum_{i=0}^{n-1}|\frac{w_i}{\sqrt{n}\lambda_{n}}|+|w_{n}|\big)\le c\big(\frac{\|w\|}{\lambda_n}+|w_n|\big).$

As a consequence, we have the following upper bound for the second term:
\begin{align}\label{II}
\mathrm{II}\le &\bigg |\bigg(\sum_{i=0}^{n}d_i w_i\bigg)a_t\bigg|\le  L c^2  (\|w\|+|w_n|)\big(\|w\|/\lambda_n+|w_n|\big)\nonumber \\ \le &   L c^2  \big(\|w\|^2/\lambda_n+2\|w\||w_n|+|w_n|^2\big).\end{align}

Finally, we proceed to estimate the third term. Since $d_n=\frac{\lambda_n^n}{\prod_{i=0 }^{n-1}(\lambda_n-\lambda_i)}>0$,  it is easy to get

\begin{align*}\mathrm{III}
=&\frac{\underline{b}-b_t}{\underline{b}}\sum_{i=0}^{n}d_i w_i\sum_{i=0}^{n}\lambda_i w_i\\
=& \frac{\underline{b}-b_t}{\underline{b}} \bigg \{\sum_{i=0}^{n-1}d_i w_i\sum_{i=0}^{n-1}\lambda_i w_i+d_{n}w_{n}\sum_{i=0}^{n-1}\lambda_i w_i
+\bigg(\sum_{i=0}^{n-1}d_i w_i\bigg)\lambda_{n} w_{n}+\lambda_n d_n w_n^2 \bigg \}
\\
\le & \frac{\underline{b}-b_t}{\underline{b}} \bigg \{\sum_{i=0}^{n-1}d_i w_i\sum_{i=0}^{n-1}\lambda_i w_i+d_{n}w_{n}\sum_{i=0}^{n-1}\lambda_i w_i
+\bigg(\sum_{i=0}^{n-1}d_i w_i\bigg)\lambda_{n} w_{n}\bigg \} .\end{align*}

Since we know that $| d_{n}|\le \frac{c}{(2n+1) \sqrt{n}}$, $|d_i|\le \frac{c}{(2n+1) n\lambda_{n}}$, $i=0,\cdots,n-1$, and $0<\lambda_i < 2n+1, i=0,\cdots,n-1$, and $\sum_{i=0}^{n-1}|w_i|\le \sqrt{n}\|w\|$,  we can easily  get the following three inequalities:
\begin{align*}
\left|\sum_{i=0}^{n-1}d_i w_i\sum_{i=0}^{n-1}\lambda_i w_i\right| \le & \sum_{i=0}^{n-1}\frac{c|w_i|}{(2n+1) n\lambda_{n}}  \sum_{i=0}^{n-1}(2n+1) |w_i|
 \\\le&  \frac{c}{ n\lambda_{n}} \bigg(\sum_{i=0}^{n-1}|w_i|\bigg)^2 \le \frac{c}{\lambda_n}\|w\|^2;
 \end{align*}

\begin{align*}
\left|d_{n}w_{n}\sum_{i=0}^{n-1}\lambda_i w_i\right|\le& \frac{c(2n+1)\sqrt{n}}{(2n+1)\sqrt{n}}|w_n| \|w\|=c|w_n|\|w\|;
\end{align*}

\begin{align*}
\bigg|\bigg(\sum_{i=0}^{n-1}d_i w_i\bigg)\lambda_n w_n\bigg| \le &\bigg(\sum_{i=0}^{n-1}\frac{c}{(2n+1) n\lambda_{n}} |w_i| \bigg)\lambda_n |w_n | \\
\le&  \frac{c}{(2n+1)\sqrt{n}}\|w\| |w_n|\le c\|w\| |w_n|.\end{align*}

Therefore, the upper bound of the  third term can be estimated as
\begin{align}\label{III}
\mathrm{III}= \frac{\underline{b}-b_t}{\underline{b}}\sum_{i=0}^{n}d_i w_i\sum_{i=0}^{n}\lambda_i w_i\nonumber
\le\frac{\overline{b}-\underline{b}}{\underline{b}}\left(\frac{c}{\lambda_{n}}\|w\|^2+2c|w_{n}|\|w\|\right).\end{align}
Denote $m=L c^2+\frac{\left(\overline{b}-\underline{b}\right)c}{\underline{b}}$.  Combining (\ref{I})-(\ref{III}), we have \begin{align}
\overset{.}{V}(\overline{w})\le& (m/\lambda_n-2)\|w\|^2+2m\|w\||w_n| -\left(\lambda_n-L c^2\right)w_n^2
\end{align}
Since $(\lambda_0,\cdots,\lambda_n) \in \Omega_{\Lambda}$, we can see $\lambda_n>\max\{2n+2, m^2+Lc^2\}$. If $m\le 1$, then $\lambda_n>2n+2>m$; if $m\geq 1$, then $\lambda_n>m^2+Lc^2>m.$ Therefore, $\lambda_n>m$ always holds whenever $(\lambda_0,\cdots,\lambda_n) \in \Omega_{\Lambda}$.

 By the inequality $\lambda_n>m$, we have $2-\frac{m}{\lambda_n}\geq 1$, therefore $$\overset{.}{V}(\overline{w}(t))\le -\|w(t)\|^2+2m\|w(t)\||w_n(t)|-(\lambda_n-L c^2)w_n^2(t).$$
Since $\lambda_{n}> m^2+L c^2$, we conclude that  $\overset{.}{V}(\overline{w}(t) )\le -\alpha \|\overline{w}(t)\|^2$ for some  $\alpha>0$, i.e., $\overset{.}{V}(\overline{w}(t) )\le -2 \alpha V(\overline{w}(t) )$.
Therefore, by the comparison theorem, we have $V(\overline{w}(t) )\le e^{-2\alpha t} V(\overline{w}(0) )$ for $t\in [0,T)$.

Finally, we  estimate the upper bound of $\left\|\overline{Y}(t)\right\|$  as follows: \begin{align*} \left\|\overline{Y}(t)\right\|& = \|P\overline{w}(t)\| \le \|P\|\|\overline{w}(t)\| = \|P\|\sqrt{2V(\overline{w}(t))} \\ \le&\|P\|\sqrt{2e^{-2\alpha t}V(\overline{w}(0))} = e^{-\alpha t}\|P\| \|\overline{w}(0)\|\\\le & e^{-\alpha t}\|P\|\|P^{-1}\|\left \|\overline{Y}(0)\right\| \le ce^{-\alpha t}\left\|\overline{Y}(0)\right\|.\end{align*}

This completes the proof of Lemma 2.$\Box$

\textbf{Proof of lemma 3.}
It suffices to show that if the parameters $(k_0,\cdots,k_n)\in \Omega_{L_0,\underline{b},b_0,c}$, then for any $T_0<T$ we have $$\left\|\overline{Y}(t)\right\|\le  c e^{-\alpha t}\left \|\overline{Y}(0)\right\|,\forall t\in [0,T_0).$$

Denote $a\overset{\triangle}{=}\sup_{0\le t\le T_0}\left\|\overline{Y}(t)\right\|$, $L'\overset{\triangle}{=}\sup_{0\le \rho \le a}\frac{\tau_2{(\rho)}}{\rho}$, $b'\overset{\triangle}{=}\tau_1(a)$.

 It is easy to verify that  $|a_t|\le \tau_2(\left\|\overline{Y}(t)\right\|) =\frac{\tau_2(\left\|\overline{Y}(t)\right\|)}{\left\|\overline{Y}(t)\right\|} \left\|\overline{Y}(t)\right\| \le L' \left\|\overline{Y}(t)\right\|$ and $\underline{b}\le b_t\le\tau_1(\left\|\overline{Y}(t)\right\|) \le \tau_1(a)=b'$ for $t\in [0,T_0)$.

Therefore, by Lemma 2, we have $\left\|\overline{Y}(t)\right\|\le  c e^{-\alpha t}\left \|\overline{Y}(0)\right\|$ for $t\in [0,T_0)$ whenever the parameters $(k_0,\cdots,k_n)\in \Omega_{L',\underline{b},b',c}$.
As a consequence, $a=\sup_{0\le t\le T_0}\left\|\overline{Y}(t)\right\|\le c\left \|\overline{Y}(0)\right\|\le c R$ and $b'=\tau_1(a)\le \tau_1(c R)$, which implies $L'\le L_0, b'\le b_0$.

The final thing we need to prove is $\Omega_{L_0,\underline{b},b_0,c} \subset \Omega_{L',\underline{b},b',c}$.

From (\ref{Omega})-(\ref{omega2}), we know that
\begin{small}
\begin{align*}\Omega_{L,\underline{b},\overline{b},c}\overset{\triangle}{=}\left\{\begin{bmatrix} k_0\\ \vdots \\k_n \end{bmatrix} \bigg|  \begin{bmatrix} k_0\\ \vdots  \\k_n\end{bmatrix}=\frac{1}{\underline{b}}\begin{bmatrix}
\prod_{i=0}^{n}\lambda_{i}\\
\vdots\\
\sum_{i=0}^{n}\lambda_{i}
\end{bmatrix} , \lambda \in \Omega_1\cap \Omega_2 \right\}
\end{align*}
\end{small}
where \begin{align}\label{168}\Omega_2=\bigg\{\lambda\in \mathbb{R}^{n+1} \bigg |  \lambda_n> \big( L c^2+\left(\overline{b}-\underline{b}\right)c/\underline{b}\big)^2+Lc^2 \bigg \}.\end{align}
 From (\ref{omega1}), it is easy to see that  $\Omega_1$ does not depend on $L,\underline{b},\overline{b}$ and $c$.
By (\ref{168}), we know that  $\Omega_2$  depends on $L,\underline{b},\overline{b},c$, i.e., $\Omega_2=\Omega_2(L,\underline{b},\overline{b},c)$. It is easy to see that if  $\underline{b}$ and $c$ are fixed, then $\Omega_2$ gets smaller for larger $L,\overline{b}$, i.e., $\Omega_2(L_0,\underline{b},b_0,c)\subset \Omega_2(L',\underline{b},b',c)$.

Therefore, we have $\Omega_{L_0,\underline{b},b_0,c} \subset \Omega_{L',\underline{b},b',c}$.
This means that, if the parameters $(k_0,\cdots,k_n)\in \Omega_{L_0,\underline{b},b_0,c}$, then for  $T_0<T$ we have $$\left\|\overline{Y}(t)\right\|\le  c e^{-\alpha t}\left \|\overline{Y}(0)\right\|,~~\forall t\in [0,T_0).$$

Since $T_0$ is arbitrary, we complete the proof of Lemma 3.$\Box$

\textbf{Proof of lemma 4.}
 Suppose that  $(k_0,\cdots,k_n)\in \Omega_{L,\underline{b}, \overline{b},c}$ and denote $\overline{Y}(t)=P\overline{w}(t) $, then similar to the proof of Lemma 2, we have \begin{equation}\label{88}
\begin{cases}
\overset{.}{w}_0&=-\lambda_0 w_0+d_0 \left( \Delta_t+ b_t\sum_{i=2}^{n}k_i \epsilon^{n-i} \xi_i\right) \\
&\vdots\\
\overset{.}{w}_{n}&=-\lambda_{n} w_{n}+d_{n} \left (\Delta_t+ b_t\sum_{i=2}^{n}k_i \epsilon^{n-i} \xi_i\right)\\
\overset{.}{\xi}&=\frac{1}{\epsilon}B \xi+ E_t
\end{cases}
\end{equation}
where $\Delta_t\overset{\triangle}{=}a_t+\frac{\underline{b}-b_t}{\underline{b}}
\sum_{i=0}^{n}\lambda_i w_i$ and $$E_t\overset{\triangle}{=}\left(0,\cdots,a_t-b_t\left(\sum_{i=0}^{n}k_i y_i-\sum_{i=2}^{n}k_i \epsilon^{n-i} \xi_i\right)\right)^T.$$

Now, we introduce a quadratic  function as follows:
$$V_0(\overline{w}(t),\xi(t))=\frac{1}{2}\|\overline{w}(t)\|^2+\xi^T(t) Q \xi(t),$$
where $Q$ is the unique positive definite matrix satisfying $B^T Q+Q B=-I$.

Then it is not difficult to compute the time derivative of $V_0$,
\begin{align}\label{49}
\overset{.}{V}_0(\overline{w}(t),\xi(t))=&
\overset{.}{V}(\overline{w}(t))+b_t\sum_{i=0}^{n}d_i w_i
\sum_{i=2}^{n}k_i \epsilon^{n-i}\xi_i
-\|\xi\|^2/\epsilon+2\xi^T Q E_t
\end{align}
where $\overset{.}{V}(\overline{w}(t))$ is computed  in (\ref{V}).

From the proof of Lemma 2, we conclude that if $(k_0,\cdots,k_n)\in \Omega_{L,\underline{b}, \overline{b},c}$, then there exists a constant $\alpha>0$, such that $\overset{.}{V}(\overline{w}(t))\le -\alpha \|\overline{w}(t)\|^2$.

Without loss of generality, assume that $\epsilon<1$, then we  have
\begin{align}\label{51}
\big|\sum_{i=2}^{n}k_i \epsilon^{n-i} \xi_i\big| &\le \sqrt{n} \max\{k_2,\cdots,k_n\} \|\xi\|,~~
\big|\sum_{i=0}^{n}\lambda_iw_i\big| &\le \sqrt{n+1}\lambda_n \|\overline{w}\|\end{align}

By  Lemma 1, we have $|d_i|\le \frac{c}{(2n+1)\sqrt{n}}, i=0,\cdots,n$, and therefore \begin{align}\label{52}\big|\sum_{i=0}^{n}d_{i}w_{i}\big|  \le c\frac{|w_0|+\cdots+|w_n|}{(2n+1)\sqrt{n}}\le c\|\overline{w}\|.\end{align}

 Hence  by $|a_t|\le L\left\|\overline{Y}(t)\right\|\le Lc \left\|\overline{w}(t)\right\|$ and $0<\underline{b}\le b_t \le \overline{b}$, we have  from (\ref{49})-(\ref{52})
 \begin{align}
 \label{50}\overset{.}{V}_0(\overline{w}(t),\xi(t)) \le& -\alpha \|\overline{w}\|^2 -\|\xi\|^2 /\epsilon +c_1\|\overline{w}\|\|\xi\|+c_2 \|\xi\|^2 \nonumber \\ =&-\alpha \|\overline{w}\|^2  +c_1\|\overline{w}\|\|\xi\|-\left(1/\epsilon-c_2\right) \|\xi\|^2 \end{align}
 where $c_1,c_2$ are two constants defined by $$c_1=\sqrt{n} \overline{b}c \max\{k_2,\cdots,k_n\} +2\lambda_{\max}(Q)\left(Lc+\overline{b}\sqrt{n+1}\lambda_n/\underline{b}\right)$$ and $c_2=2\sqrt{n}\overline{b} \lambda_{\max}(Q)\max\{k_2,\cdots,k_n\}$, which are independent of $\epsilon$.

 Since $\alpha$, $c_1,c_2$ are  independent of $\epsilon$, it is easy to see   from (\ref{50}) that there exists $\epsilon^*>0$ such that  whenever $0<\epsilon<\epsilon^*$, we have
 $$\overset{.}{V}_0(\overline{w}(t),\xi(t))\le -2\beta V_0(\overline{w}(t),\xi(t)) $$ for some $\beta>0.$
By the comparison theorem, we have  $$V_0(\overline{w}(t),\xi(t))\le e^{-2\beta t} V_0(\overline{w}(0),\xi(0)).$$ As a consequence, we have
\begin{align*}
\|Y(t)\|\le & \|\overline{Y}(t)\| =\|P\overline{w}(t)\| \le \|P\|\|\overline{w}(t)\|\\
\le& \|P\|\sqrt{2V_0(\overline{w}(t),\xi(t))} \le \|P\|e^{-\beta t}\sqrt{2V_0(\overline{w}(0),\xi(0))}\\
 =&\|P\|e^{-\beta t}\left(\|\overline{w}(0)\|+\sqrt{2\lambda_{\max}(Q)}\|\xi(0)\|\right) \\
\le &\|P\| e^{-\beta t} (\left\|P^{-1}\right\|\left\|\overline{Y}(0)\right\|+\sqrt{2\lambda_{\max}(Q)}\|\xi(0)\|)\\
 \le &ce^{-\beta t}(\|Y(0)\|+|y_0(0)|+\sqrt{2\lambda_{\max}(Q)}\|\xi(0)\|).
 \end{align*}
Moreover, we have
\begin{align*}
&\sqrt{\lambda_{\min}(Q)}\|\xi(t)\|  \le \sqrt{V_0(\overline{w}(t),\xi(t))}
\le e^{-\beta t} \sqrt{V_0(\overline{w}(0),\xi(0))}\\
\le& e^{-\beta t} \left( \left\|\overline{w}(0)\right\|+ \sqrt{\lambda_{\max}(Q)}\|\xi(0)\|\right)
\le  e^{-\beta t}\left(\left\|P^{-1}\right\|\left\|\overline{Y}(0)\right\|+\sqrt{\lambda_{\max}(Q)}\|\xi(0)\|\right)\\
 \le& e^{-\beta t}\left(c\left\|\overline{Y}(0)\right\|+\sqrt{\lambda_{\max}(Q)}\|\xi(0)\|\right)\Box\end{align*}

\subsection{Appendix B}

\textbf{Proof of Example 2.}
 Obviously, (\ref{triangular}) has uniform relative degree $2$. By simple calculations, we have  $\Phi(x)=(x_1,f_1(x_1)+x_2)^T$, $F(x)=f_1'(x_1)(f_1(x_1)+x_2)+f_2(x_1,x_2)$ and $G(x)=1$ for any $x\in\mathbb{R}^2$, which implies that Assumption (A) is satisfied.

 Next,   the Jacobian matrix $J_\Phi(x)=\begin{bmatrix}1&0\\f_1'(x_1)&1\end{bmatrix}$ is nonsingular and it is easy to see that $\lim_{\|x\|\to \infty}\|\Phi(x)\|=\infty$. Therefore, by Theorem A1, we conclude that  Assumption (C') is satisfied.

Let $y^*$ be any given setpoint. Then  $x^*\overset{\triangle}{=}\Phi^{-1}(y^*,0)=(y^*,-f_1(y^*))$. Now we proceed to estimate the upper bound of $|F(x)G(x^*)-F(x^*)G(x)|$,
  \begin{align*}&|F(x)G(x^*)-F(x^*)G(x)|\\
  =&|F(x)-F(x^*)|=|F(x_1,x_2)-F(y^*,-f_1(y^*))|\\
  = &|f_1'(x_1)(f_1(x_1)+x_2)+f_2(x_1,x_2)-f_2(y^*,-f_1(y^*))|\\
  \le &|f_1'(x_1)(f_1(x_1)+x_2)|+|f_2(x_1,x_2)-f_2(y^*,x_2)|+|f_2(y^*,x_2)-f_2(y^*,-f_1(y^*))|\\
  \le &L |f_1(x_1)+x_2|+L|x_1-y^*|+L|x_2+f_1(y^*)|\\
  \le &L |f_1(x_1)+x_2|+L|x_1-y^*|+L(|f_1(x_1)+x_2|+|f_1(y^*)-f_1(x_1)|)\\
   \le &2L |f_1(x_1)+x_2|+ (L+L^2)|x_1-y^*|\\
   \le &\sqrt{4L^2+(L+L^2)^2} \sqrt{(f_1(x_1)+x_2)^2+(x_1-y^*)^2}\\
   = &\sqrt{4L^2+(L+L^2)^2} \|\Phi(x)-\Phi(x^*)\|
  \end{align*}
Therefore, by Theorem  2 and Remark 5, the classical PID controller $u(t)=k_p e(t)+k_i \int_{0}^{t}e(s)ds +k_d \frac{de(t)}{dt}$  can globally stabilize the system and make the regulation error $e(t)$ converge to $0$ exponentially for any initial state $x(0)\in\mathbb{R}^2$ as long as $(k_i,k_p,k_d)\in \Omega_{\tilde{L},1,1,c}$, where $\tilde{L}=\sqrt{4L^2+(L+L^2)^2}$.$\Box$

\textbf{Proof of Example 3.}
 By some simple calculations, it is easy to get $\Phi(x)=(x_1,f_1(x_1,x_2))$, $F(x)=\frac{\partial f_1}{\partial x_1}(x)f_1(x)+\frac{\partial f_1}{\partial x_2}(x)f_2(x)$, $G(x)=\frac{\partial f_1}{\partial x_2}(x)g(x)\geq \underline{b}b_1>0$ and therefore Assumption (A) is satisfied.

 Let  $y^*$ be any given setpoint. Since $\frac{\partial f_1}{\partial x_2}(x)\geq \underline{b}>0$ for any $x\in\mathbb{R}^2$, then there exists a unique  $x_2^*\in\mathbb{R}$ such that  $f_1(y^*,x_2^*)=0$. Denote $x^*\overset{\triangle}{=}\Phi^{-1}(y^*,0)=(y^*,x_2^*)$, $\overline{x}_1\overset{\triangle}{=}x_1-y^*$ and $\overline{x}_2\overset{\triangle}{=}x_2-x_2^*$. Then we can easily  obtain the following equalities:
\begin{align}\label{x*}
\|x-x^*\|^2 =(x_1-y^*)^2+(x_2-x_2^*)^2=\overline{x}_1^2+\overline{x}_2^2
\end{align}and
$\|\Phi(x)-\Phi(x^*)\|^2= \overline{x}_1^2+(f_1(x_1,x_2)-f_1(y^*,x_2^*))^2 .$
By the mean value theorem, we have
\begin{align}\label{z*}\|\Phi(x)-\Phi(x^*)\|^2= \overline{x}_1^2+ (\theta_1 \overline{x}_1+\theta_2 \overline{x}_2)^2,\end{align}
where $|\theta_1|\le L$, $0<\underline{b}\le \theta_2\le L$.

Now we proceed to prove  $\|H(x)-H(x^*)\|\le L_0 \|\Phi(x)-z^*\|$ for some $L_0>0.$
First, we will show that \begin{align}\label{aa}\alpha \|x-x^*\|\le \|\Phi(x)-\Phi(x^*)\|, \forall x \in \mathbb{R}^2\end{align} for some $\alpha>0$.

Without loss of generality, assume that $x-x^*\neq 0$. Then there exists some $r>0$ and $\theta\in [0,2\pi]$ such that $\overline{x}_1=r \cos{\theta}$ and $\overline{x}_2=r \sin{\theta}$ . From (\ref{x*})-(\ref{z*}), we have
\begin{align*}
\|\Phi(x)-\Phi(x^*)\|^2/ \|x-x^*\|^2=\cos^2\theta+(\theta_1 \cos \theta+\theta_2\sin \theta)^2.
\end{align*}
Denote $\alpha=\inf\sqrt{\cos^2\theta+(\theta_1 \cos \theta+\theta_2\sin \theta)^2}$, where the infimum is taken for all  $\theta\in [0,2\pi],|\theta_1|\le L, 0<\underline{b}\le \theta_2\le L$.  It is easy to obtain that  $\alpha>0$, i.e., (\ref{aa}) is satisfied.

Next, we will prove that \begin{align}\label{bb}\|H(x)-H(x^*)\|\le \beta \|x-x^*\|\end{align} for some   $\beta>0$.

We only give a sketch proof due to space limitation. First, note that $f_1(x^*)=0$ and $|f_1(0)|\le M$, it is not difficult to get the upper bound of $\|x^*\|$ by the assumption $\left\|\frac{\partial f_1}{\partial x}\right\|\le L$ and $\frac{\partial f_1}{\partial x_2}(x)\geq \underline{b}>0$. Then by $|f_2(0)|\le M$ and the upper bound of $\|x^*\|$, it is not difficult to  find $M_0$ such that $|f_2(x^*)|\le M_0.$

 Recall that $f_1(x^*)=0$, then we have
\begin{align*}&|F(x)-F(x^*)|
\\=&\left|\frac{\partial f_1}{\partial x_1}(x)f_1(x)\right|+\left|\frac{\partial f_1}{\partial x_2}(x)f_2(x)-\frac{\partial f_1}{\partial x_2}(x^*)f_2(x^*)\right|
\\\le& \left|\frac{\partial f_1}{\partial x_1}(x)f_1(x)\right|+\left|\frac{\partial f_1}{\partial x_2}(x)(f_2(x)-f_2(x^*))\right|
+\left|\left(\frac{\partial f_1}{\partial x_2}(x)-\frac{\partial f_1}{\partial x_2}(x^*)\right)f_2(x^*)\right|
\\\le &L^2 \|x-x^*\|+\left|\frac{\partial f_1}{\partial x_2}(x)(f_2(x)-f_2(x^*))\right|
+\left|\left(\frac{\partial f_1}{\partial x_2}(x)-\frac{\partial f_1}{\partial x_2}(x^*)\right)f_2(x^*)\right|
\\\le &\left(2L^2+M_0L\right)\|x-x^*\|\end{align*}

Similarly, we can obtain
\begin{align*}|G(x)-G(x^*)|=&\left|\frac{\partial f_1}{\partial x_2}(x)g(x)-\frac{\partial f_1}{\partial x_2}(x^*)g(x^*)\right|\\=&\left|\left(\frac{\partial f_1}{\partial x_2}(x)-\frac{\partial f_1}{\partial x_2}(x^*)\right)g(x)\right|+\left|\frac{\partial f_1}{\partial x_2}(x^*)(g(x)-g(x^*))\right|
\\\le& \left(L^2+b_2L\right)\|x-x^*\|.\end{align*}
Therefore, we have \begin{align}\label{L0}\|H(x)-H(x^*)\|\le L_0 \|\Phi(x)-\Phi(x^*)\|\end{align} for some $L_0>0.$
Furthermore, we have
\begin{align}\label{F1}
 |F(x)|\le& |F(x^*)|+|F(x)-F(x^*)|
\le \big|\frac{\partial f_1}{\partial x_2}(x^*)f_2(x^*)\big|+\|H(x)-H(x^*)\| \nonumber\\
\le& LM_0+L_0\|\Phi(x)-\Phi(x^*)\|
\le  L_0\|\Phi(x)\|+LM_0+L_0|y^*|.\end{align}
By (\ref{L0})-(\ref{F1}) and $\underline{b}b_1\le G(x)=\frac{\partial f_1}{\partial x_2}(x)g(x)\le L b_2$, we conclude that Assumption (B1) is satisfied.

Finally, by the fact $\frac{\partial f_1}{\partial x_2}(x)\geq \underline{b}>0$, it is easy to see that $\lim_{\|x\|\to \infty}\|\Phi(x)\|=\infty$.  By Theorem A1, we conclude that Assumption (C') is also satisfied. Therefore, by Theorem  2, the classical PID controller $u(t)=k_p e(t)+k_i \int_{0}^{t}e(s)ds +k_d \frac{de(t)}{dt}$  can globally stabilize the system and make the regulation error $e(t)$ converge to $0$ exponentially for any initial state $x(0)\in\mathbb{R}^2.$ $\Box$

\textbf{Proof of Remark 4.}  We will give an example to show that the  super-linear growth rate (\ref{13}) in Assumption (C) cannot be weakened to (\ref{exa}) for any $\eta>0$.

We first define a function $f$  as follows:
\begin{align}f(x)=
\begin{cases}
2-\left(\log x\right)^{-\eta}, ~x\geq e\\
-f(-x)~~,~~~ x\le -e\\
\end{cases},
\end{align} where $e$ is the natural logarithm. We can extend $f$ as a smooth function defined on $\mathbb{R}$ with $f(0)=0$ and $f'(x)> 0,$ when $-e<x<e$. Let us consider the following nonlinear uncertain plant with PID controller:
\begin{align}\label{eta}
\begin{cases}
\overset{.}{x}_1&=\epsilon f(x_2)\\
\overset{.}{x}_2&=\frac{1+u}{\epsilon f'(x_2)}\\
~ y &= x_1\\
u(t)&=k_p e(t)+k_i \int_{0}^{t}  e(s)ds+k_d\frac{de}{dt}(t)
\end{cases}
\end{align} where $\epsilon$ is an unknown constant with $0<\epsilon \le1.$
Let $y^*=0$ be the setpoint. Then  we can verify that both Assumptions (A) and  (B) are satisfied and  (\ref{exa}) is satisfied for some positive numbers $N_1$ and $N_2$. However, we will show that for any $R>0$ and for any given  PID parameters,  there always exists initial state  $x(0)\in\mathbb{R}^2$ satisfying $\|x(0)\|\le R$, such that the solution of the closed-loop system (\ref{eta}) with the initial state $x(0)$ has  finite escape time for all $\epsilon<\min\{\frac{1}{4(|k_p|+|k_i|+|k_d|)},\frac{1}{4}\}$.

First, by some simple calculations,  we can obtain $\Phi(x)=(x_1,\epsilon f(x_2))$, $F(x)=1,G(x)=1$ for any $x\in\mathbb{R}^2$.

Define an increasing function $\tau_1(r)=r+2$. It is not difficult to see that $\sup_{x_2 \in\mathbb{R}}|f(x_2)|=2$.

It is easy to obtain  $$\|\Phi(x)\|\le |x_1|+|\epsilon f(x_2)|\le \|x\|+2 =\tau_1(\|x\|)$$ and $$\|H(x)\|\le \sqrt{2}<\tau_1(\|x\|)$$ for any $x\in \mathbb{R}^2.$ For $y^*=0$, then $x^*=(0,0)$. We have $$\|H(x)-H(x^*)\|=0\le \tau_2(r)=r.$$
 Therefore, Assumptions (A) and (B) are satisfied.

By simple calculations, we can obtain $$J_\Phi^{-1}(x)=\begin{bmatrix}1&0\\0&\frac{1}{\epsilon f'(x_2)}\end{bmatrix}.$$
By the definition of $f$, it is easy to see that $$f'(x_2)=\eta\left(|x_2|\log^{1+\eta} |x_2|\right)^{-1}, \forall |x_2|\geq e$$
Therefore, we can obtain
$$\big\|J_\Phi^{-1}(x)\big\| \le N_1\|x\|\log^{1+\eta}\|x\|+N_2,\forall x\in\mathbb{R}^2$$
for some constants $N_1$ and $N_2$(possibly depend on $\epsilon$).

Let $R>0$ and the parameter triple $(k_p,k_i,k_d)\in\mathbb{R}^3$ are given arbitrarily. Let $[0,T)$ be the maximal existence interval of (\ref{eta}) with initial state $x(0)=(0,0)$.  We proceed to prove that the closed-loop equation (\ref{eta}) will have finite escape time for $\epsilon$ sufficiently small, i.e., $T<\infty$.

Denote $T_0=\min\{T,1\}$. Let us first prove that $u(t)=O(\epsilon)$ for $t\in[0,T_0)$.

Now suppose that $t\in[0,T_0)$, then
\begin{align}\label{ee}|e(t)|=&|x_1(t)|=\left|x_1(0)+\int_{0}^t \overset{.}{x}_1(s)ds\right|=\left|\int_{0}^t \epsilon f(x_2(s))ds\right|\le \int_{0}^t 2\epsilon ds\le 2\epsilon.\end{align}
By (\ref{ee}), we have
\begin{align}\label{eee}\left|\int_{0}^t e(s)ds\right|\le \int_{0}^t |e(s)|ds\le 2\epsilon t\le  2\epsilon \end{align}
Recall that $|f(x)|\le 2$, we can obtain
\begin{align}\label{ed}\left|\overset{.}{e}(t)\right|=\left|\overset{.}{x}_1(t)\right|=|\epsilon f(x_2(t))|\le 2\epsilon. \end{align}
From (\ref{ee})-(\ref{ed}), we have  \begin{align}\label{pid}|u(t)|=&\left|k_p e(t)+k_i \int_{0}^{t}  e(s)ds+k_d\frac{de}{dt}(t)\right|\le 2(|k_p|+|k_i|+|k_d|)\epsilon.\end{align}
Next, we proceed to prove that $T<\infty$ whenever $\epsilon<\min\{\frac{1}{4(|k_p|+|k_i|+|k_d|)},\frac{1}{4}\}$.

From $\overset{.}{x}_2(t)=\frac{1+u(t)}{\epsilon f'(x_2(t))}$ and (\ref{pid}), we conclude that if $\epsilon<\frac{1}{4(|k_p|+|k_i|+|k_d|)}$, then we have $\overset{.}{x}_2\geq \frac{1}{2\epsilon f'(x_2)}$ for $t\in [0,T_0),$ i.e., $f'(x_2) d{x}_2\geq \frac{1}{2\epsilon }dt$, $~\forall t\in [0,T_0)$.

By the comparison lemma in differential equations, the following inequality will be satisfied \begin{align}\label{fx}f(x_2(t))-f(x_2(0))=f(x_2(t))\geq \frac{t}{2\epsilon},~~ t\in [0,T_0).\end{align}
Notice that $|f(x)|\le 2$, therefore by (\ref{fx}), we know that  $T_0 \le 4\epsilon$.

Since $\epsilon<\frac{1}{4}$, we obtain $T_0 \le 4\epsilon<1$, which implies $T<1<\infty$, i.e.,  the maximal existence interval $[0,T)$ of the closed-loop system (\ref{eta}) with initial state $(0,0)$ is finite.$\Box$

\bibliographystyle{plain}

\end{document}